\documentclass[10pt]{article}
\usepackage{amsfonts}
\usepackage{amsmath,amssymb}
\usepackage[dvips]{graphicx}
\usepackage{algorithm}
\usepackage{algorithmicx}
\usepackage{algpseudocode}
\usepackage{yhmath}
\usepackage{qtree}
\usepackage{xcolor}
\usepackage{epsfig}
\usepackage{lineno,hyperref}
\usepackage{mathtools}
\usepackage{amsmath,bm}
\usepackage{amsthm}
\allowdisplaybreaks
\usepackage[T1]{fontenc}
\usepackage[utf8]{inputenc}
\usepackage{moresize}
\usepackage{makecell}

\usepackage{booktabs}

\usepackage{authblk}
\usepackage{bigints}

\usepackage[toc]{appendix}
\numberwithin{equation}{section}

\usepackage{verbatim}
\usepackage{graphicx}
\usepackage{geometry}
\usepackage{epsfig}
\usepackage{lineno,hyperref}
\usepackage[many]{tcolorbox}
\usetikzlibrary{shadows}
\usepackage{caption}
\newtcolorbox{shadedbox}{
  breakable,
  enhanced jigsaw,
  colback=white,
}

\newcommand{\be}{\begin{equation}}
\newcommand{\ee}{\end{equation}}
\newcommand{\ba}{\begin{array}}
\newcommand{\ea}{\end{array}}
\newcommand{\bea}{\begin{eqnarray*}}
\newcommand{\eea}{\end{eqnarray*}}
\newcommand{\bean}{\begin{eqnarray}}
\newcommand{\eean}{\end{eqnarray}}

\newtheorem{theorem}{Theorem}[section]
\newtheorem{lemma}{Lemma}[section]
\newtheorem{remark}{Remark}[section]

\newtheorem{corollary}{Corollary}[section]
\newtheorem{example}{Example}[section]

\newcommand\numberthis[1][]{%
    \refstepcounter{equation}%
    \ifx#1\empty\else\label{eq:#1}\fi%
    \tag{\theequation}%
}

\makeatletter
\def\BState{\State\hskip-\ALG@thistlm}
\makeatother

\newcommand{\lc}{\mathrel{\raise2pt\hbox{${\mathop<\limits_{\raise1pt\hbox{\mbox{$\sim$}}}}$}}}
\newcommand{\gc}{\mathrel{\raise2pt\hbox{${\mathop>\limits_{\raise1pt\hbox{\mbox{$\sim$}}}}$}}}
\newcommand{\ec}{\mathrel{\raise1pt\hbox{${\mathop=\limits_{\raise2pt\hbox{\mbox{$\sim$}}}}$}}}

\renewcommand{\theequation}{\arabic{section}.\arabic{equation}}

\usepackage{fancyhdr}
\begin{document}

\title{Energy Stable L2 Schemes for Time-Fractional Phase-Field Equations}
\author[1]{Chaoyu Quan}
%\author[2,1]{Tao Tang}
\author[2]{Boyi Wang}
%\author[3,4]{Jiang Yang}
\affil[1]{\small SUSTech International Center for Mathematics, Southern University of Science and Technology, Shenzhen, China (\href{mailto:quancy@sustech.edu.cn}{quancy@sustech.edu.cn}).}
%\affil[2]{\small Division of Science and Technology, BNU-HKBU United International College, Zhuhai, Guangdong, China (\href{mailto:tangt@sustech.edu.cn}{tangt@sustech.edu.cn}).}
\affil[2]{\small Department of Mathematics, National University of Singapore, Singapore 119076 \& Department of Mathematics, Southern University of Science and Technology, Shenzhen 518055, China(\href{mailto:e0220126@u.nus.edu}{e0220126@u.nus.edu}).}
%\affil[4]{\small Guangdong Provincial Key Laboratory of Computational Science and Material Design, Southern University of Science and Technology, Shenzhen, China.}
\maketitle

\begin{abstract}
In this article, the energy stability of two high-order L2 schemes for time-fractional phase-field equations is established. 
We propose a reformulation of the L2 operator and also some new properties on it. 
We prove the energy boundedness (by initial energy) of an L2 scalar auxiliary variable scheme for any phase-field equation and the fractional energy law of an implicit-explicit L2 Adams--Bashforth scheme for the Allen--Cahn equation. 
The stability analysis is based on a new Cholesky decomposition proposed recently by some of us. 
%Some numerical experiments are provided to verify the  stability analysis.
\end{abstract}


{\bf AMS: } 65M06, 65M12, 74A50

\section{Introduction}
Phase-field models have been widely-used in many areas, such as material sciences, multiphase flows, biology, and image processing, etc. 
One important feature of phase-field model is that its energy admits a dissipation law with respect to time.
In particular, this property has become a criterion for designing numerical schemes for phase-field equations in the past decade.

From the numerical point of view, the resolution of phase-field equation is interesting and challenging due to the existence of nonlinearity. 
Moreover, it is usually expected that the maximum principle and the energy dissipation could be preserved for a numerical scheme of phase-field equation. 
So far, there have been different energy stable schemes including the convex-splitting scheme \cite{eyre1998unconditionally,wang2010unconditionally}, the stabilization scheme \cite{shen2010numerical,xu2006stability}, and the scalar auxiliary variable (SAV) scheme \cite{shen2018scalar}.

In this article, we study the energy dissipation property of high order schemes for phase-field models with Caputo time-derivative.
The time-fractional phase-field equation can be written in the general form of
\begin{equation}\label{eq:phase_field}
	 \partial^\alpha_t  u =  \mathcal G \mu,
\end{equation}
where $\alpha \in (0,1)$, $\mathcal G$ is a nonpositive operator depending on the phase-field model, $\mu = \delta_u E$ is the functional derivative of some energy $E$, and $\partial_t^\alpha$ is the Caputo derivative \cite{caputo1967linear} defined by
\begin{equation}\label{eq:fracder}
\partial^\alpha_t  u(t) \coloneqq \frac{1}{\Gamma(1-\alpha)}\int_0^t \frac{ u'(s)}{(t-s)^\alpha}\, {\rm d} s, \qquad t\in(0,T),
\end{equation}
with $\Gamma(\cdot)$ the gamma function.
Taking different functional $\mathcal G$ and $\mu$, \eqref{eq:phase_field} becomes different phase-field equation, such as the Allen--Cahn (AC) model \cite{allen1979microscopic}, the Cahn--Hilliard (CH) model \cite{cahn1958free} and the molecular beam epitaxy (MBE) model \cite{li2003thin}.
%More specifically, in the AC model \cite{allen1979microscopic} and the CH  model \cite{cahn1958free}, $\mu$ is taken to be
%\begin{equation}\label{eq:acch}
%\mu = -\varepsilon^2 \Delta   u + F'( u),
%\end{equation}
%where $\varepsilon>0$ is the interface width parameter and $F$ is a double-well potential functional.
%We take the common choice $F( u)  = \frac 1 4 \left(1- u^2\right)^2$ in this article so that $ F'( u) =  u^3- u$.
For the sake of simplicity, we consider  the periodic boundary condition for the time-fractional phase-field equation \eqref{eq:phase_field}.

Straightforward computation of the derivative of energy with respect to time gives
\begin{equation}\label{eq:ed}
\frac{\rm d}{{\rm d} t} E( u)
=  \int_\Omega \partial_t  u  \left({\mathcal G^{-1}}\partial^\alpha_t  u  \right) {\rm d}x.
\end{equation}
It is known that when $\alpha=1$, i.e., the conventional case,  the phase-field models are gradient flows. 
So the energy associated with these models decays with time, that is the so-called energy dissipation law. 
However, it is still unknown if such energy dissipation property holds in the general case of $0<\alpha<1$.

In \cite{TangYZ19}, the authors demonstrated that the classical energy of \eqref{eq:phase_field} is bounded from above by the initial energy. 
Later, it is observed numerically in \cite{du2020time} and then proved theoretically in \cite{CSIAM-AM-1-478} that the time-fractional derivative of energy is always nonpositive, i.e., the so-called fractional energy law, 
\begin{equation}
\partial_t^\alpha E(t)\leq 0,\quad \forall\, 0<t<T.
\end{equation}
Moreover, discrete fractional energy law has been obtained in \cite{quan2020numerical} for first and $2-\alpha$ order schemes.
For example, for first-order L1 schemes, the discrete fractional energy law is satisfied 
\begin{equation}
\overline\partial_n^\alpha E  \coloneqq \sum_{k=1}^{n} b_{n-k} D_k E \leq 0\quad \forall n\geq 1,
\end{equation}
where
\begin{equation}\label{eq:Dju}
b_j = \frac{\Delta t^{1-\alpha}}{ \Gamma(2-\alpha)}\left[ (j+1)^{1-\alpha} - j^{1-\alpha} \right]
\quad \mbox{and}\quad
D_j u \coloneqq \frac{u^{j}-u^{j-1}}{\Delta t},
\quad j\geq 0.
\end{equation}
See for example \cite{sun2006fully,lin2007finite} for the deviation and analysis of L1 coefficients $b_j$. 
In addition, there are some other interesting works on time-fractional gradient flows. 
For example, Li and Salgado develop the theory of fractional gradient flows that minimize a convex l.s.c. energy in \cite{li2021time}; Fritz, Khristenko, and Wohlmuth propose the equivalence between a time-fractional and a integer-order gradient flow in 
 \cite{fritz2021equivalence} where a dissipation-preserving augmented energy is introduced.

It is natural to generalize the energy stability analysis to higher-order schemes.
In this work, we consider two L2 schemes \cite{lv2016error}: one is a second order L2 SAV scheme for any phase field equation and the other is a $3-\alpha$ order implicit-explicit L2 Adams--Bashforth (AB) scheme for the Allen--Cahn equation. 
We prove that the energy of the L2 SAV scheme for any phase-field equation is bounded by initial energy. 
Moreover, the implicit-explicit L2 AB scheme satisfied the fractional energy law, i.e., the fractional derivative of energy is nonpositive. 
In fact, the analysis is based on two new properties of the L2 operator $L_k^\alpha$:
\begin{equation}
\sum_{k = 1}^n \left< L_{k}^\alpha u, 3 D_k u - D_{k-1} u\right> \geq 0,
\end{equation}
and
\begin{equation}
\sum_{k = 1}^n  d_{n-k+1} \left< L_{k}^\alpha u, D_k u\right>  \geq 0,
\end{equation}
where the definitions of $L_k^\alpha$ and $d_j$ are given in Section \ref{sect2}.

This article is organized as follows.
In Section \ref{sect2}, we propose a reformulation of L2 approximation and then prove the aforementioned properties of L2 operator. 
In Section \ref{sect3}, we study the energy stability of  an implicit-explicit L2 AB scheme and an L2 SAV scheme. 
Some numerical tests are given in Section \ref{sect4}.
Finally, we give a brief conclusion in the last section.
%The widely-used L1 approximation \cite{lin2007finite} of time fractional derivative \eqref{eq:fracder} at $t_k$ is 
%\begin{equation}\label{eq:partial_t_alpha}
%D_k^\alpha  u \coloneqq \sum_{j=1}^{k} b_{k-j} D_{j}  u,
%\end{equation}
%where $D_j$ is the discrete first-order derivative
%\begin{equation}\label{eq:Dju}
%D_j u \coloneqq \frac{u^{j}-u^{j-1}}{\Delta t}
%\quad \mbox{and}\quad
%b_j = \frac{\Delta t^{1-\alpha}}{ \Gamma(2-\alpha)}\left[ (j+1)^{1-\alpha} - j^{1-\alpha} \right],
%\quad j\geq 0.
%\end{equation}
%As stated in \cite{quan2020numerical}, the L1 approximation of fractional derivative \eqref{eq:partial_t_alpha} satisfies,
%\begin{equation}
%\sum_{k=1}^{n} b_{n-k}  \left<D_k^\alpha u, D_k u\right> \geq 0 \quad \forall n\geq 1.
%\end{equation}
%This yields the discrete fractional energy law, i.e.,
%\begin{equation}
%D_n^\alpha E  = \sum_{k=1}^{n} b_{n-k} D_k E \leq 0\quad \forall n\geq 1,
%\end{equation}
%for the L1 schemes of time-fractional phase-field equations as soon as $D_k E\leq \left<\mu^k, D_k u\right>$ always holds.

%The above energy law for the $2-\alpha$ order scheme \eqref{eq:L1SAV} is essentially based on the property of L1 operator $D_{n+\frac 1 2}^\alpha$ that is similar to the property of $D_{n}^\alpha$ studied in \cite{quan2020numerical}.

\section{Analysis of L2 approximation}\label{sect2}
In this section, we prove some useful properties of the L2 operator $ L_n^\alpha$. 

Let $\Delta t = T/N$ be the time step size and $t_k = k \Delta t$, $0\leq k \leq N$.
The L2 approximation \cite{lv2016error} of time fractional derivative \eqref{eq:fracder} is written as 
\begin{equation}\label{eq:L2}
  \begin{array}{r@{}l}
	\begin{aligned}
	 L_1^\alpha  u  & = &&\frac{1}{\Gamma(2-\alpha)\Delta t^{\alpha}} \left( u^1- u^0\right),\quad k = 1,\\
 	 L_k^\alpha  u & = &&\frac{1}{\Gamma(3-\alpha)\Delta t^\alpha} {\Bigg \{} \sum_{j=1}^{k-1} \left(a_j u^{k-j-1}+b_j u^{k-j}+c_j u^{k-j+1}\right) \\ 
	 & &&+  \frac{\alpha}{2} u^{k-2}-2 u^{k-1}+\frac{4-\alpha}{2}u^k {\Bigg\}},\quad k\geq 2,\\	
	\end{aligned}
  \end{array}
\end{equation}
where 
\begin{equation}\label{eq:abc}
  \begin{array}{r@{}l}
	\begin{aligned}
	a_j & = -\frac{3}{2}(2-\alpha)(j+1)^{1-\alpha}+\frac 1 2 (2-\alpha)j^{1-\alpha} + (j+1)^{2-\alpha}-j^{2-\alpha}, \\
	b_j & = 2(2-\alpha)(j+1)^{1-\alpha}-2(j+1)^{2-\alpha}+2j^{2-\alpha},\\
	c_j & = -\frac 1 2 (2-\alpha)\left((j+1)^{1-\alpha}+j^{1-\alpha}\right)+(j+1)^{2-\alpha}-j^{2-\alpha}.
	\end{aligned}
  \end{array}
\end{equation}
Note that the relationship $a_j+b_j+c_j = 0$ holds.

\subsection{Reformulation of L2 operator}
Why shall we reformulate the L2 coefficients in \eqref{eq:L2}?
The reason is that $b_j$ is not monotonic w.r.t. $j$, which leads to the difficulty when analyzing the positive-definiteness property of L2 operator. 

We propose to reformulate \eqref{eq:L2} as 
\begin{equation}\label{eq:L2reform}
  \begin{array}{r@{}l}
	\begin{aligned}
	 L_1^\alpha  u & = \frac{\Delta t^{1-\alpha}}{\Gamma(3-\alpha)} {\Big \{} r_1 D_1 u + d_1 D_1 u{\Big\}}, \quad k = 1,\\  
 	 L_k^\alpha  u %& = \frac{\Delta t^{1-\alpha}}{\Gamma(3-\alpha)} {\bigg \{} \frac{4-\alpha}{2}D_k u - \frac{\alpha}{2} D_{k-1} u+ \sum_{j=1}^{k-1} \left(c_j D_{k-j+1} u - a_j D_{k-j} u \right){\bigg\}}\\
	 & = \frac{\Delta t^{1-\alpha}}{\Gamma(3-\alpha)} \bigg\{ \frac{3\alpha}{2} D_k u - \frac{\alpha}{2} D_{k-1} u    +  \sum_{j=1}^{k} d_j D_{k-j+1} u  -c_kD_{1} u \bigg\}, \quad k \geq 2,
	\end{aligned}
  \end{array}
\end{equation}
where $D_j u$ is defined in \eqref{eq:Dju},  
\begin{equation}\label{eq:r1}
r_1 = 2-\alpha - d_1 = 2 + \frac 1 2 \alpha - \left(\frac \alpha 2 +1\right) 2^{1-\alpha}> \frac 3 4 \alpha,\qquad \alpha\in (0,1),
\end{equation}
and
\begin{equation}\label{eq:dj}
d_j = \left\{
  \begin{array}{r@{}l}
	\begin{aligned}
	& c_1 + 2-2\alpha, &&  j = 1,\\
	& c_j - a_{j-1}, && j = 2,\ldots,k. \\
	\end{aligned}
  \end{array}	
  \right.
\end{equation}
To be precise, we can write $d_j$ as
\begin{equation}\label{eq:djform0}
  \begin{array}{r@{}l}
	\begin{aligned}
	d_1 & = \left(1+\frac\alpha2\right) 2^{1-\alpha} -\frac32\alpha, \quad j =1,\\
	d_j & = \left(1-\frac\alpha 2\right) \left[ -(j+1)^{1-\alpha}+2j^{1-\alpha}-(j-1)^{1-\alpha}\right] + \left[ (j+1)^{2-\alpha}-2j^{2-\alpha}+(j-1)^{2-\alpha} \right] \\
	& =  - \left(1-\frac\alpha 2\right)  \kappa(j,1-\alpha) + \kappa(j,2-\alpha),\quad j\geq 2,
	\end{aligned}
  \end{array}
\end{equation}
where 
\begin{equation}\label{eq:kappa0}
\kappa(j,\beta) \coloneqq (j+1)^{\beta}-2j^{\beta}+(j-1)^{\beta}.
\end{equation}

Now we propose the following properties of $a_j,~c_j,$ and $d_j$ that will be useful in our later energy analysis.
\begin{lemma}[Properties of L2 operator]\label{lem1}
For any $\alpha\in(0,1)$, the following properties on the L2 coefficients $a_j,~c_j,~d_j$ hold:
\begin{itemize}
\item[(1)] $a_j<0$, $a_j-a_{j+1}<0$, and $3a_j-4a_{j+1}+a_{j+2}<0$ increase w.r.t $j$;
\item[(2)]$c_j>0$, $c_j-c_{j+1}>0$, and $3c_j-4c_{j+1}+c_{j+2}>0$ decrease w.r.t. $j$; 
\item[(3)] $d_j>0$, $d_j -d_{j+1} >0$, and $3d_j-4d_{j+1}+d_{j+2}>0$ decrease w.r.t. $j$;
\item[(4)] $4 d_{j+1}\geq d_j$.
\end{itemize}
\end{lemma}
\begin{proof}
We prove the above properties one by one.
We treat the index $j\geq 1$ as a continuous variable so that the derivatives w.r.t. $j$ can be computed.

(1) From \cite[Eq. (2.3)]{lv2016error} and variable transformation, $a_j$ can be written in the integral form of 
\begin{equation}
a_j = \frac{(2-\alpha)(1-\alpha) \Delta t^\alpha}{2\Delta t^2} \int_0^{\Delta t} \frac{2s-3\Delta t}{(j\Delta t +\Delta t-s)^\alpha} \,{\rm d}s< 0. 
\end{equation}
It is easy to find that 
\begin{equation}
\partial_j a_j >0 \quad\mbox{and}\quad \partial_{jj} a_j <0,
\end{equation}
implying $a_j<0$ and $a_j-a_{j+1}<0$ increases.
Furthermore, we have
\begin{equation}
3a_j-4a_{j+1}+a_{j+2}  = \frac{(2-\alpha)(1-\alpha) \Delta t^\alpha}{2\Delta t^2} \int_0^{\Delta t} (2s-3\Delta t)\rho(j,s)\,{\rm d}s< 0
\end{equation}
with 
\begin{equation}\label{eq:rhojs}
\rho(j,s) = 3(j\Delta t +\Delta t-s)^{-\alpha} -4((j+1)\Delta t +\Delta t-s)^{-\alpha} +((j+2)\Delta t +\Delta t-s)^{-\alpha}.
\end{equation}
It is not difficult to verify $\rho(j,s)>0$ and $\partial_j\rho(j,s)< 0$, which yields that
\begin{equation}
\partial_j \left(3a_j-4a_{j+1}+a_{j+2}\right) > 0. 
\end{equation}

(2) Similarly, $c_j$ can be written in the integral form of 
\begin{equation}
c_j = \frac{(2-\alpha)(1-\alpha) \Delta t^\alpha}{2\Delta t^2} \int_0^{\Delta t} \frac{2s-\Delta t}{(j\Delta t +\Delta t-s)^\alpha} \,{\rm d}s >0.
\end{equation}
Then we have
\begin{equation}
\partial_j c_j <0 \quad\mbox{and}\quad \partial_{jj} c_j >0,
\end{equation}
implying $c_j>0$ and $c_j-c_{j+1}>0$ decreases.
Furthermore, we have
\begin{equation}
3c_j-4c_{j+1}+c_{j+2}  = \frac{(2-\alpha)(1-\alpha) \Delta t^\alpha}{2\Delta t^2} \int_0^{\Delta t} (2s-\Delta t)\rho(j,s)\,{\rm d}s< 0
\end{equation}
with $\rho(j,s)>0$ given by \eqref{eq:rhojs} satisfying $\partial_j\rho(j,s)< 0$ and $\partial_s\rho(j,s) > 0$.
As a consequence, we have
\begin{equation}
\partial_j \left(3c_j-4c_{j+1}+c_{j+2}\right) < 0. 
\end{equation}

(3) According to the above properties of $a_j$ and $c_j$, $d_j >0$, $d_j -d_{j+1}>0$, and $3d_j -4d_{j+1}+d_{j+2}>0$  decrease w.r.t. $j$ when $j\geq 2$.
Moreover, when $j=1$, straight computation gives
\begin{equation}
\begin{aligned}
& d_1 - d_2 = (c_1-c_2) + a_1+2-2\alpha> 0, \\
& d_1 - d_2 - (d_2 - d_3) = c_1-2c_2+c_3 +2a_1-a_2 +2-2\alpha>0,\\
& 3 d_1 - 4d_2+d_3 - (3d_2 -4 d_3+d_4) > 0.
\end{aligned}
\end{equation}

(4) In the case of $j=1$, we can obtain 
\begin{equation}
4 d_2-d_1 = 4(c_2-a_1)-c_1-2+2\alpha = 2(4+\alpha)3^{1-\alpha}-\frac 9 2\left(2+\alpha\right) 2^{1-\alpha} + \frac 7 2\alpha > 0.
\end{equation}
In the case $2\leq j\leq n-1$, 
\begin{equation}\label{eq:djform}
  \begin{array}{r@{}l}
	\begin{aligned}
	d_j 
	& =  - \left(1-\frac\alpha 2\right)  \kappa(j,1-\alpha) + \kappa(j,2-\alpha),
	\end{aligned}
  \end{array}
\end{equation}
where 
\begin{equation}\label{eq:kappa}
\kappa(j,\beta) \coloneqq (j+1)^{\beta}-2j^{\beta}+(j-1)^{\beta}.
\end{equation}
Due to the concavity of $j^{1-\alpha}$ and the convexity of $j^{2-\alpha}$, it is easy to see
\begin{equation}
\kappa(j,1-\alpha) < 0 \quad\mbox{and}\quad \kappa(j,2-\alpha)>0.
\end{equation}
According to the Jensen's inequality, the following inequality holds
\begin{equation}\label{eq:rho1}
  \begin{array}{r@{}l}
	\begin{aligned}
 -4\kappa(j+1,1-\alpha)+\kappa(j,1-\alpha) 
& = -4(j+2)^{1-\alpha}+9(j+1)^{1-\alpha} - 6j^{1-\alpha}+(j-1)^{1-\alpha} \\
& \geq -j^{1-\alpha}+(j-1)^{1-\alpha}.
	\end{aligned}
  \end{array}
\end{equation}
Similarly, we also have
\begin{equation}\label{eq:rho2}
4\kappa(j+1,2-\alpha)-\kappa(j,2-\alpha) \geq j^{2-\alpha}-(j-1)^{2-\alpha}\geq j^{1-\alpha}-(j-1)^{1-\alpha}.
\end{equation}
Combining \eqref{eq:djform}, \eqref{eq:rho1}, and \eqref{eq:rho2}, we obtain
\begin{equation}
4 d_{j+1}-d_{j} \geq 0,\quad \forall 2\leq j\leq n-1.
\end{equation}
%In the case of $j = n-1$, we have
%\begin{equation}
%  \begin{array}{r@{}l}
%	\begin{aligned}
%	4d_{n}-d_{n-1} 
%	& = 4(-a_{n-1}) -(c_{n-1}-a_{n-2}) \\
%	& = \frac 1 2 (2-\alpha) \left[ 13 n^{1-\alpha} - 6 (n-1)^{1-\alpha} + (n-2)^{1-\alpha} \right] \\
%	& + \left[ -5n^{2-\alpha}+6(n-1)^{2-\alpha}-(n-2)^{2-\alpha}\right]\\
%	& {\red \geq 0.}
%	\end{aligned}
%  \end{array}
%\end{equation}
In summary, we conclude that $4 d_{j+1}\geq d_j$, $\forall 1\leq j\leq n-1$. 
\end{proof}

\subsection{Positive definiteness}
Based on Lemma \ref{lem1}, we first state and prove the following theorem on the discrete operator $L_t^\alpha$ given by \eqref{eq:L2reform}. 

\begin{lemma}\label{main-lem}
% u \in C^1([0,T]; H^2(\Omega)\cap H^1_0(\Omega))
For any function $u\in C\left([0,T]; L_2(\Omega)\right)$, the following inequality on the operator $L_k^\alpha$ holds:
\begin{equation}\label{eq:lem2.2}
\sum_{k = 1}^n \left< L_{k}^\alpha u, 3 D_k u - D_{k-1} u\right>\geq \frac {\alpha \Delta t^{1-\alpha}} {2\Gamma(3-\alpha)} \sum_{k=1}^{n} \left\| D_k u\right\|^2 \geq 0.
\end{equation}
\end{lemma}
\begin{proof}
According to the formula \eqref{eq:L2} of $L^\alpha_k u$, we can write the left-hand side of \eqref{eq:lem2.2} in the following matrix form:
\begin{equation}\label{eq:discrete_ED}
  \begin{array}{r@{}l}
	\begin{aligned}
	\sum_{k = 1}^n \left< L_{k}^\alpha u, 3 D_k u - D_{k-1} u\right>
	& =  \frac{\Delta t^{1-\alpha}}{\Gamma(3-\alpha)} \int_{\Omega}\psi \left(\mathbf A+\mathbf B +\mathbf C\right)\psi^{\rm T}\, {\rm d} x,
 	\end{aligned}
  \end{array}
\end{equation}
with 
\begin{equation}
\begin{aligned}
&\psi = \left[D_1 u,D_2 u,\cdots,D_n u\right],\\
& \mathbf A =
\left[\begin{array}{cccccc} 
\frac 1 2 (3d_1 + a_2) &  &  &  & &  
\\ -3a_2+a_3 & \frac 1 2 (3 d_1 - d_2) & & &  &  
\\ -3a_3+a_4 & 3d_2-d_3 & \frac 1 2 (3 d_1 - d_2)  & &  &  
\\\vdots & \vdots & \vdots & \ddots &   & 
\\ -3 a_{n-1}+a_{n} & 3 d_{n-2}-d_{n-1} & 3 d_{n-3}-d_{n-2} & \cdots & \frac 1 2 (3 d_1 - d_2) &  
\\ -3 a_n & 3 d_{n-1} & 3d_{n-2} &\cdots & 3d_2 & \frac 5 2 d_1\end{array}\right],\\
&\mathbf B =
\left[\begin{array}{ccccc} \frac 1 2 (3d_1+a_2)  &  &  &  &  \\ -d_1 & \frac 1 2 (3d_1-d_2)  &  &  &  \\  & \ddots & \ddots &  &  \\   &  & -d_1 & \frac 1 2 (3d_1-d_2) &  \\  &   &  & -d_1 & \frac 1 2 d_1  \end{array}\right],\\
%\left[\begin{array}{cccccc}  d_1 &  &  &  &  & \\ -d_1 &  d_1 &  &  &  & \\ 
% & \ddots & \ddots &  &  &\\ 
% &  & \ddots &  \ddots &  &\\ 
%    & &  & -d_1 & d_1  &  \\
%      -a_{n+1} & d_n  & \cdots & d_4 & d_3-d_1 & d_1+d_2  \end{array}\right],
&\mathbf C =
\left[\begin{array}{ccccc} r_1 &  &  &  &  \\ -\frac 1 2 \alpha & \frac 3 2 \alpha  &  &  &  \\  & \ddots & \ddots &  &  \\   &  & -\frac 1 2 \alpha & \frac 3 2 \alpha &  \\  &   &  & -\frac 1 2 \alpha & \frac 3 2 \alpha  \end{array}\right].
\end{aligned}
\end{equation}
Here we make a split of the associated matrix which will facilitate the proof. 

On the right-hand side of \eqref{eq:discrete_ED}, we actually split the essential matrix into three matrices $\mathbf A$, $\mathbf B$, and $\mathbf C$. 
It is not difficult to see that $\mathbf B$ is positive definite since $\frac 1 2 (3d_1 + a_2)>d_1$ and $\frac 1 2 (3d_1 -d_2) >d_1$, and $\mathbf C$ is also positive definite due to $r_1 > \frac 3 4 \alpha$. 
Further, $\mathbf C$ satisfies
\begin{equation}\label{ineq:C}
\psi \mathbf C \psi^{\rm T} \geq \frac\alpha 2 \psi \psi^{\rm T}.
\end{equation}
As a consequence, to derive \eqref{eq:lem2.2}, the remaining work is to prove that $\mathbf A$ is definite positive, which is equivalent to prove that $\mathbf M = \mathbf A + \mathbf A^{\rm T}$ is positive definite.

To prove the positive definiteness of $\mathbf M$, we split it into
\begin{equation}
  \begin{array}{r@{}l}
	\begin{aligned}
\mathbf M = \mathbf A + \mathbf A^{\rm T}
& = \left[\begin{array}{cc} \mathbf M_{n-1} & \mathbf b^{\rm T} \\ \mathbf b & 5d_1\end{array}\right],
%& =
%\left[\begin{array}{cccccc} 
%3 d_1 + a_2 & -3a_2+a_3 & -3a_3+a_4  & \cdots & -3 a_{n-1}+a_{n}  &  -3 a_n 
%\\ -3a_2+a_3 & 3 d_1 - d_2 & 3d_2-d_3 & \cdots & 3 d_{n-2}-d_{n-1}  &  3 d_{n-1} 
%\\ -3a_3+a_4 & 3d_2-d_3 & 3 d_1 - d_2  & \cdots & 3 d_{n-3}-d_{n-2} &  3d_{n-2}
%\\\vdots & \vdots & \vdots & \ddots &   & 
%\\ -3 a_{n-1}+a_{n} & 3 d_{n-2}-d_{n-1} & 3 d_{n-3}-d_{n-2} & \cdots & 3 d_1 - d_2 &  3d_2 
%\\ -3 a_n & 3 d_{n-1} & 3d_{n-2} &\cdots & 3d_2 & 5 d_1\end{array}\right]. \\
 	\end{aligned}
  \end{array}
\end{equation}
where $\mathbf M_{n-1}$ is the leading principle minor of $\mathbf M$ of size $(n-1)\times (n-1)$. 
Note that $0<-a_j<d_j$ holds true and $\mathbf M$ is a symmetric matrix composed of positive elements.
According to Lemma \ref{lem1}, $\mathbf M_{n-1}$ satisfies the three conditions in \cite[Lemma 2.1]{CSIAM-AM-1-478}: for the lower triangular part of $\mathbf M_{n-1}$,
\begin{equation}
\begin{aligned}
& ({\rm i})~ [\mathbf M_{n-1}]_{i-1,j}\geq [\mathbf M_{n-1}]_{i, j}; \\
& ({\rm ii})~[\mathbf M_{n-1}]_{i, j-1}< [\mathbf M_{n-1}]_{i, j}; \\
& ({\rm iii})~ [\mathbf M_{n-1}]_{i-1, j-1} - [\mathbf M_{n-1}]_{i, j-1}\leq [\mathbf M_{n-1}]_{i-1, j} - [\mathbf M_{n-1}]_{i, j}.
\end{aligned}
\end{equation}
Therefore it has a Cholesky decomposition
\begin{equation}
\mathbf M_{n-1} = \mathbf L_{n-1}\mathbf L_{n-1}^{\rm T} ,
\end{equation}
where the lower triangular part of $\mathbf L_{n-1}$ is composed of positive elements decreasing along each column.
Further, based on Lemma \ref{lem1}, we can find the following matrix 
\begin{equation}
 \begin{array}{r@{}l}
	\begin{aligned}
\widetilde{\mathbf M}
& = \left[\begin{array}{cc} \mathbf M_{n-1} & \frac 2 3 \mathbf b^{\rm T} \\ \frac 2 3 \mathbf b & 2 d_1\end{array}\right] 
 	\end{aligned}
  \end{array}
\end{equation}
also satisfies the three conditions in \cite[Lemma 2.1]{CSIAM-AM-1-478} and can be decomposed as
\begin{equation}
\widetilde {\mathbf M} = 
\left[\begin{array}{cc} \mathbf L_{n-1} & \\ {\mathbf l} & {l}_{nn}\end{array}\right]
\left[\begin{array}{cc} \mathbf L_{n-1}^{\rm T} & {\mathbf l}^{\rm T}  \\  & { l}_{nn}\end{array}\right],
\end{equation}
where the lower triangular matrix on the right-hand side satisfies the properties in \cite[Lemma 2.1]{CSIAM-AM-1-478}.
The following inequality holds:
\begin{equation}
{\mathbf l}\,{\mathbf l}^{\rm T} = 2d_1- l_{nn}^2 < 2d_1.
\end{equation}
Therefore, we can derive
\begin{equation}
 \begin{array}{r@{}l}
	\begin{aligned}
{\mathbf M}
& = \left[\begin{array}{cc} \mathbf M_{n-1} & \mathbf b^{\rm T} \\  \mathbf b & 5 d_1\end{array}\right] 
=
\left[\begin{array}{cc} \mathbf L_{n-1} & \\ \frac 3 2 {\mathbf l} & {l}_{nn}\end{array}\right]
\left[\begin{array}{cc} \mathbf L_{n-1}^{\rm T} & \frac 3 2 {\mathbf l}^{\rm T}  \\  & {l}_{nn}\end{array}\right] .
 	\end{aligned}
  \end{array}
\end{equation}
Note that
\begin{equation}
l_{nn}^2 = 5d_1 - \frac 9 4 {\mathbf l}\,{\mathbf l}^{\rm T} > 5d_1 - \frac 9 2 d_1 = \frac 1 2 d_1 >0.
\end{equation}
This implies that the above decomposition is feasible and one can take $l_{nn} > 0$.
We have proven that $\mathbf M$ is positive definite and so is $\mathbf A$. 
In summary, $\mathbf A,~\mathbf B,$ and $\mathbf C $ are all positive definite. 
Combining \eqref{eq:discrete_ED} and \eqref{ineq:C}, we then have \eqref{eq:lem2.2}.
The proof is completed.
\end{proof}

Furthermore, we state and prove the following theorem on the discrete operator $L_t^\alpha$ given by \eqref{eq:L2reform}. 

\begin{lemma}\label{main-lem2}
% u \in C^1([0,T]; H^2(\Omega)\cap H^1_0(\Omega))
For any function $u\in C\left([0,T]; L_2(\Omega)\right)$, the following inequality on the operator $L_k^\alpha$ holds:
\begin{equation}\label{eq:thm3.1}
\sum_{k = 1}^n  d_{n-k+1} \left< L_{k}^\alpha u, D_k u\right> \geq \frac {5\alpha \Delta t^{1-\alpha}} {12\Gamma(3-\alpha)} \sum_{k=1}^{n}  d_{n-k+1} \left\| D_k u\right\|^2 \ \geq 0.
\end{equation}
\end{lemma}
\begin{proof}
According to the formula of $L^\alpha_t u^k$, we have
\begin{equation}\label{eq:discrete_ED_2}
  \begin{array}{r@{}l}
	\begin{aligned}
	 \sum_{k = 1}^n  d_{n-k+1} \left< L_k^\alpha u, D_k u\right> 
	& =  \frac{\Delta t^{1-\alpha}}{\Gamma(3-\alpha)} \int_{\Omega}\psi \left(\mathbf A+\mathbf B\right)\psi^{\rm T}\, {\rm d} x,
 	\end{aligned}
  \end{array}
\end{equation}
with 
\begin{equation}
\psi = \left[D_1 u,D_2 u,\cdots,D_n u\right],
\end{equation}
\begin{equation}
\mathbf A =
\left[\begin{array}{ccccc} d_n &  &  &  &  \\ & d_{n-1}  &  &  &  \\ &  & \ddots &  &  \\ &  &  & d_2 &  \\ &  &  &  & d_1 \end{array}\right]
\left[\begin{array}{ccccc} d_1 &  &  &  &  \\ -a_1 & d_1 &  &  &  \\\vdots & \vdots & \ddots &  &  \\ -a_{n-2} & d_{n-2} & \cdots & d_1 &  \\ -a_{n-1} & d_{n-1} & \cdots & d_2 & d_1\end{array}\right],
\end{equation}
and
\begin{equation}
\mathbf B =
\left[\begin{array}{ccccc} d_n &  &  &  &  \\ & d_{n-1} &  &  &  \\ &  & \ddots &  &  \\ &  &  & d_2 &  \\ &  &  &  & d_1 \end{array}\right]
\left[\begin{array}{ccccc} r_1 &  &  &  &  \\ -\frac 1 2 \alpha & \frac 3 2 \alpha  &  &  &  \\  & \ddots & \ddots &  &  \\   &  & -\frac 1 2 \alpha & \frac 3 2 \alpha &  \\  &   &  & -\frac 1 2 \alpha & \frac 3 2 \alpha  \end{array}\right].
\end{equation}
%It is sufficient to prove that $\mathbf A$ and $\mathbf B$ are both definite positive matrices.

We first prove that $\mathbf B$ is strictly positive definite.
It is not difficult to verify that $r_1> \frac 3 4 \alpha$ as pointed out in \eqref{eq:r1}.
In Lemma \ref{lem1}, we have proven that $d_{j} \geq \frac 1  4 d_{j-1}$.
As a consequence, we have
\begin{equation}\label{ineq:B}
  \begin{array}{r@{}l}
	\begin{aligned}
	\psi \mathbf B \psi^{\rm T} 
	& = r_1 d_n  \psi_1^2 + \sum_{j = 2}^n \left(\frac{3\alpha}{2}  d_{n-j+1}  \psi_j^2 -\frac \alpha 2 d_{n-j+1} \psi_{j-1} \psi_j\right) \\
	& \geq \frac {3\alpha} 4  d_n \psi_1^2 + \alpha \sum_{j=2}^n  \left(\frac 3 2 d_{n-j+1}  \psi_j^2 -\frac 1 2 d_{n-j+1} \psi_{j-1} \psi_j\right) \\
	& \geq \frac {5\alpha}{12} d_n \psi_1^2  +   \alpha \sum_{j=2}^n  d_{n-j+1}  \left(\frac 1 {12}  \psi_{j-1}^2 + \frac 7 6 \psi_j^2- \frac 1 2 \psi_{j-1} \psi_j\right) \\
	& = \frac {5\alpha}{12} d_n \psi_1^2  +   \alpha \sum_{j=2}^n  d_{n-j+1}  \left(\frac 1 {12}  (\psi_{j-1}-3\psi_j)^2 + \frac 5 {12} \psi_j^2 \right) \\
	& \geq \frac {5\alpha} {12} \sum_{j=1}^{n}  d_{n-j+1ß} \psi_{j}^2.
 	\end{aligned}
  \end{array}
\end{equation} 

Next, we prove that $\mathbf A$ is positive definite, which is equivalent to prove that $\mathbf A + \mathbf A^{\rm T}$ is positive definite.
We consider the following conjugate transformation of $\mathbf A + \mathbf A^{\rm T}$:
\begin{equation}\label{eq:conj}
 {\mathbf S} = P \left( {\mathbf A} +{\mathbf A}^{\rm T} \right) P^{\rm T},
\end{equation}
where $P$ is an anti-diagonal matrix
\begin{equation}\label{eq:Pn}
P = \left[ \begin{array}{cccc} &  &  & d_{1}^{-1} \\ &  & d_{2}^{-1} &  \\ & \adots &  &  \\ d_{n}^{-1} &  &  & \end{array}\right]_{n\times n}.
\end{equation}
As a consequence, the lower triangular part of $\mathbf S$ can be written in the form of
\begin{equation}
\mathbf S_{ij}  = \left\{
  \begin{array}{r@{}l}
	\begin{aligned}
	& 2 d_{1}d_{i}^{-1} && \mbox{if } i = j,\\
	& d_{i-j+1}d_{i}^{-1}  && \mbox{if } j<i<n,\\
	& -a_{n-j}d_{n}^{-1}  && \mbox{if } j<i=n.
	 \end{aligned}
  \end{array}
  \right.
\end{equation}
Note that $0<-a_i<d_i$ holds true and $\mathbf S$ is a symmetric matrix composed of positive elements.
%We split $\mathbf M$ into two symmetric matrix
%\begin{equation}
%\mathbf M = \mathbf S + \mathbf T
%\end{equation}
%where $\mathbf S$ and $\mathbf T$ are symmetric with lower triangular parts: $i\geq j$,
%\begin{equation}
%\mathbf S_{ij}  = \left\{
%  \begin{array}{r@{}l}
%	\begin{aligned}
%	& 2d_{1}d_{i}^{-1} && \mbox{if } i = j,\\
%	& d_{i-j+1}d_{i}^{-1}  && \mbox{if } i > j,
%	 \end{aligned}
%  \end{array}
%  \right.
%\end{equation}
%and
%\begin{equation}
%\mathbf T_{ij}  = \left\{
%  \begin{array}{r@{}l}
%	\begin{aligned}
%	& c_i d_{i}^{-1}  && \mbox{if } i > j=1,\\
%	& 0  && \mbox{else.}
%	 \end{aligned}
%  \end{array}
%  \right.
%\end{equation}

We show that the lower triangular part of $\mathbf S$ satisfies the following properties: 
\begin{equation}\label{eq:Mpro}
\begin{aligned}
& ({\rm i})~ \mathbf S_{i-1,j}\geq \mathbf S_{i, j}; \\
&({\rm ii})~\mathbf S_{i, j-1}< \mathbf S_{i, j}; \\
& ({\rm iii})~ \mathbf S_{i-1, j-1} - \mathbf S_{i, j-1}\leq \mathbf S_{i-1, j} - \mathbf S_{i, j}.
\end{aligned}
\end{equation}
From Lemma \ref{lem1}, it is easy to see that if $i\geq j$, $\mathbf S_{ij}$ increases w.r.t. $j$.
The second property in \eqref{eq:Mpro} is satisfied.
In the following proof, we treat $i$ and $j$ as variable. 
We want to prove that for all $i> j\geq 1$,
\begin{equation}\label{eq:piM}
  \begin{array}{r@{}l}
	\begin{aligned}
	\partial_i \left(d_{i-j+1}d_{i}^{-1} \right) = d_i^{-2}\left(d_i\,\partial_i d_{i-j+1}- d_{i-j+1}\,\partial_i d_i\right)\leq 0,
	 \end{aligned}
  \end{array}
\end{equation}
and
\begin{equation}\label{eq:pijM}
  \begin{array}{r@{}l}
	\begin{aligned}
	\partial_{ij} \left(d_{i-j+1}d_{i}^{-1} \right)  = d_i^{-2}\left(-d_i\,\partial_{ii} d_{i-j+1}+ \partial_i d_i \, \partial_i d_{i-j+1}\right) \leq 0.
	 \end{aligned}
  \end{array}
\end{equation}
When $j = 1$, it is clear that $\partial_i \left(d_{i-j+1}d_{i}^{-1} \right)  = 0$, which indicates that \eqref{eq:pijM} can lead to \eqref{eq:piM}.
So, we only need to prove \eqref{eq:pijM}.
Note that
\begin{equation}\label{eq:pdform}
  \begin{array}{r@{}l}
	\begin{aligned}
	d_i & = -\frac 1 2\left(2-\alpha\right) \kappa(i,1-\alpha) + \kappa(i,2-\alpha), \\
	\partial_i d_i
	 & = -\frac 1 2 (2-\alpha)(1-\alpha) \kappa(i,-\alpha) + (2-\alpha)\kappa(i,1-\alpha), \\
	 \partial_{ii} d_i 
	&= \frac \alpha 2 (2-\alpha)(1-\alpha)\kappa(i,-\alpha-1)+ (2-\alpha)(1-\alpha) \kappa(i,-\alpha),
	 \end{aligned}
  \end{array}
\end{equation}
where $\kappa(\cdot,\cdot)$ is given by \eqref{eq:kappa}.
We then have 
\begin{equation}
  \begin{array}{r@{}l}
	\begin{aligned}
	& -d_i\,\partial_{ii} d_{i-j+1}+ \partial_i d_i \, \partial_i d_{i-j+1} \\
	& =(2-\alpha)(1-\alpha) \left[\frac 1 2\left(2-\alpha \right) \kappa(i,1-\alpha) - \kappa(i,2-\alpha)\right]
	 \left[\frac \alpha 2\kappa(i-j+1,-\alpha-1)+ \kappa(i-j+1,-\alpha)\right]\\
	 & + (2-\alpha)^2 \left[\frac 1 2(1-\alpha) \kappa(i,-\alpha) - \kappa(i,1-\alpha)\right]
	 \left[\frac 1 2(1-\alpha) \kappa(i-j+1,-\alpha) - \kappa(i-j+1,1-\alpha)\right] \\
	 & = -\frac 1 2 (2-\alpha)^2(1-\alpha) Q_1
	 + \frac 1 2 (2-\alpha)^2(1-\alpha) Q_2
	 + (2-\alpha) Q_3,
	 \end{aligned}
  \end{array}
\end{equation}
where 
\begin{equation}\label{eq:Q}
  \begin{array}{r@{}l}
	\begin{aligned}
	Q_1 & =  \frac \alpha 2 \kappa(i,2-\alpha) \kappa(i-j+1,-\alpha-1) + \kappa(i,-\alpha) \kappa(i-j+1,1-\alpha),\\
	Q_2 & = \frac \alpha 2 \kappa(i,1-\alpha) \kappa(i-j+1,-\alpha-1)+ \frac 1 2 (1-\alpha) \kappa(i,-\alpha) \kappa(i-j+1,-\alpha),\\
	Q_3 & = -(1-\alpha) \kappa(i,2-\alpha) \kappa(i-j+1,-\alpha) + (2-\alpha) \kappa(i,1-\alpha) \kappa(i-j+1,1-\alpha).
	 \end{aligned}
  \end{array}
\end{equation}
In Appendix \ref{append}, we prove that $Q_1\geq 0$, $Q_2\leq 0$, and $Q_3\leq 0$, which is very technical (see Figure \ref{fig:Q} for numerical verification).
Now we can say that \eqref{eq:piM} and \eqref{eq:pijM} holds true, which implies that the three properties \eqref{eq:Mpro} are satisfied when $i<n$.
When $i = n$, using the fact that $c_{n-j}d_n^{-1}$ increases w.r.t. $j$ as well as \eqref{eq:piM} and \eqref{eq:pijM}, one can verify that the three properties \eqref{eq:Mpro} are still satisfied. 
Therefore, $\mathbf S$ is positive definite.

\begin{figure}[!h]
\centering
\includegraphics[trim ={1.5in 0 1.8in 0},clip,width=1.\textwidth]{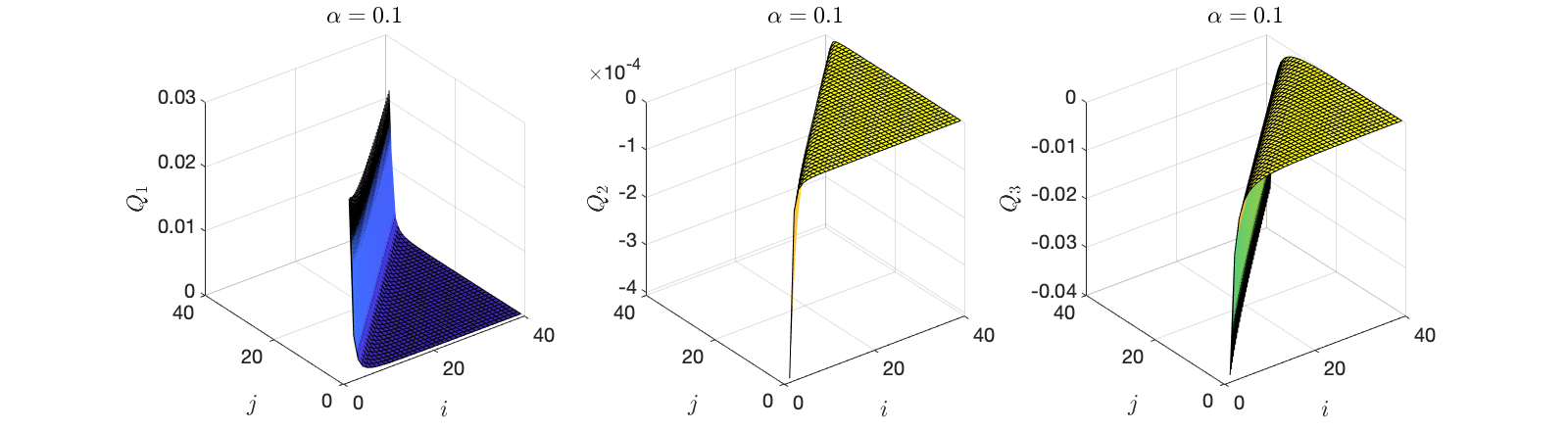}
\includegraphics[trim ={1.5in 0 1.8in 0},clip,width=1.\textwidth]{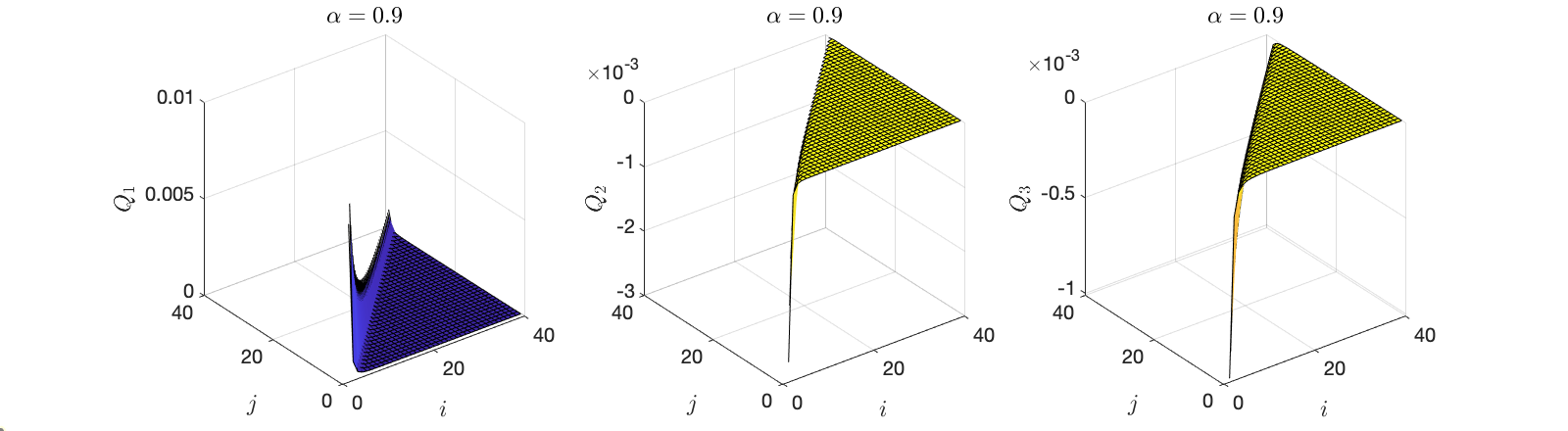}
\caption{Signs of $Q_1,~Q_2,~Q_3$ for $\alpha=0.1$ and $0.9$.}\label{fig:Q}
\end{figure}

In summary,  $\mathbf A$ is also positive definite.
As a consequence, by combining \eqref{eq:discrete_ED_2} and \eqref{ineq:B}, the inequality \eqref{eq:thm3.1} is true. 
\end{proof}

\section{Energy stable L2 schemes}\label{sect3}
In this section, we propose second order and $3-\alpha$ order schemes for time-fractional phase-field equations and establish the corresponding energy stability based on the analysis of L2 operators.

\subsection{L2 SAV scheme}
We propose a second order semi-discrete scheme for the , using the L2 approximation for the fractional derivative and the SAV technique \cite{shen2019new} for the nonlinear term:
\begin{equation}\label{eq:L2SAV}
  \begin{array}{r@{}l}
	\begin{aligned}
L_{n} ^\alpha u& =  \mathcal G \left[\mathcal L u^{n} + \frac{r^{n}}{\sqrt{E_1(\overline u^{n})}} \delta_u E_1(\overline u^{n})\right], \\
3r^{n}-4r^{n-1}+r^{n-2} & = \frac{1}{2\sqrt{E_1(\overline u^{n})}} \left<\delta_u E_1(\overline u^{n}), 3u^{n}-4u^{n-1}+u^{n-2}\right>,
	\end{aligned}
  \end{array}
\end{equation}
with $\overline u^n = 2u^{n-1}-u^{n-2}$.
Then, we can state the energy boundedness for the scheme \eqref{eq:L2SAV}.

\begin{theorem}[Energy boundedness]
For the second order L2 scheme \eqref{eq:L2SAV}, the following energy boundedness holds: $\forall 1\leq n \leq N,$
\begin{equation}
\widetilde E^{n} \leq \widetilde E^0,  
\end{equation}
where 
\begin{equation}\label{eq:modenergy}
\widetilde E^{n} = \frac 1 4 \left(\left< u^n, \mathcal L u^n\right> + \left< 2u^n-u^{n-1}, \mathcal L (2u^n-u^{n-1})\right> \right)+ \frac 1 2 \left( (r^n)^2 + (2r^n-r^{n-1})^2\right).
\end{equation}
\end{theorem}
\begin{proof}
Take the inner products of the first two equations in \eqref{eq:L2SAV} respectively with $3u^{n}-4u^{n-1}+u^{n-2}$ and $r^n$.
Then, multiply the third equation in \eqref{eq:L2SAV} with $2 r^n $.
Combining the derived three equations, we have
\begin{equation}
\left<\mathcal G^{-1}L_n^\alpha u,3u^{n}-4u^{n-1}+u^{n-2}\right> = \left<\mathcal L u^n, 3u^{n}-4u^{n-1}+u^{n-2}\right> + 2 r^n \left(3r^{n}-4r^{n-1}-r^{n-2}\right).
\end{equation}
As a consequence, we can derive
\begin{equation}
\widetilde E^n -\widetilde E^{n-1} \leq \frac 1 {2} \left<\mathcal G^{-1}L_n^\alpha u,3u^{n}-4u^{n-1}+u^{n-2}\right>.
\end{equation}
According to Lemma \ref{main-lem}, we then have
\begin{equation}
\widetilde E^n -\widetilde E^0 \leq  \frac 1 {2} \sum_{k= 1}^n  \left<\mathcal G^{-1}L_k^\alpha u,3u^{k}-4u^{k-1}+u^{k-2}\right> \leq 0.
\end{equation}
\end{proof}

\subsection{$3-\alpha$ order implicit-explicit L2 scheme}
We consider the following $3-\alpha$ order implicit-explicit L2 scheme for the time-fractional Allen--Cahn equation with $\mathcal G = -1,~\mathcal L = -\varepsilon^2 \Delta$:
\begin{equation}\label{eq:L2AB}
  \begin{array}{r@{}l}
	\begin{aligned}
L_{n+1} ^\alpha u& =\varepsilon^2 \Delta u^{n+1} - 3 f(u^n)+3f(u^{n-1})-f(u^{n-2}),
	\end{aligned}
  \end{array}
\end{equation}
where  $f(u) = u^3-u$.
%If not, a truncation technique can be used, see for example \cite{shen2010numerical}.
Then, we state the following fractional energy law for scheme \eqref{eq:L2AB} under a mild restriction on $\Delta t$.
\begin{theorem}[Fractional energy law]\label{main-thm}
For the numerical scheme \eqref{eq:L2AB}, assume that there exists a constant $L_0\geq1$ s.t.
\begin{equation}
\|u^n\|_\infty \leq L_0,\qquad\forall  n\geq 1.
\end{equation}
If
\begin{equation}
\Delta t^\alpha\leq \frac{5\alpha}{168\Gamma(3-\alpha)(3L_0-1)},
\end{equation}
then the following time-fractional energy law holds for all $n$:
\begin{equation}\label{main-thm-eq}
\sum_{k = 1}^{n}  d_{n-k+1} D_k E\leq 0,
\end{equation}
where $d_j>0$ is given by \eqref{eq:dj}.
\end{theorem}
\begin{proof}  
%Define $D_k E = \left(E^{k}+\frac L2\|u^{k}-u^{k-1}\|^2\right)-(E^{k-1}+\frac L2\|u^{k-1}-u^{k-2}\|^2)$. The following inequality will be shown
%\begin{align}
%	(E^{k}+\frac L2\|u^{k}-u^{k-1}\|^2)-(E^{k-1}+\frac L2\|u^{k-1}-u^{k-2}\|^2) \leq -\langle L_{k}^\alpha u,&u^{k}-u^{k-1}\rangle.
%\end{align}
%for all $n\geq 1$. When $k=1$, we have
%	\begin{align}\label{E-induction1}
%		E^{1}+\frac L2\|u^{1}-u^0\|^2-E^0 \leq 0.
%	\end{align}
Rewrite \eqref{eq:L2AB} as 
\begin{align}
	L_{k} ^\alpha u& =\varepsilon^2 \Delta u^{k} - 3 f(u^{k-1})+3 f(u^{k-2})- f(u^{k-3}), \quad\forall k = 1,\cdots, n.
\end{align} 
Multiplying equation  by $u^{k}-u^{k-1}$ and integrating the resultant equation over $\Omega$. We compute each term in the equation as follows
\begin{equation}
\begin{aligned}
	&\langle\varepsilon^2\Delta u^{k},u^{k}-u^{k-1}\rangle = -\frac{\varepsilon^2}{2}\left(\|\nabla u^{k}\|^2-\|\nabla u^{k-1}\|^2+\|\nabla u^{k}-\nabla u^{k-1}\|^2\right), \\
	&-\langle f(u^{k-1}),u^{k}-u^{k-1}\rangle = -\langle F(u^{k})-F(u^{k-1}),1\rangle+\frac12\langle f'(\xi_1)(u^{k}-u^{k-1}),u^{k}-u^{k-1}\rangle,\\
	&-2\langle f(u^{k-1})-f(u^{k-2}),u^{k}-u^{k-1}\rangle = -2\langle f'(\xi_2)(u^{k-1}-u^{k-2}),u^{k}-u^{k-1}\rangle,\\
	& \langle f(u^{k-2})-f(u^{k-3}),u^{k}-u^{k-1}\rangle = \langle f'(\xi_3)(u^{k-2}-u^{k-3}),u^{k}-u^{k-1}\rangle.
\end{aligned}
\end{equation}
where $\xi_i$ is between $u^{k-i}$ and $u^{k-i+1}$, $i=1,2,3$.
Summing up all equations and using $|f'(\xi_i)|\leq L = 3L_0-1$, we arrive at
\begin{align}
&\langle L_{k}^\alpha u,  u^k-u^{k-1}\rangle\leq -(E^{k}-E^{k-1})+2L\|u^k-u^{k-1} \|^2+ L \|u^{k-1}-u^{k-2}\|^2+\frac L2 \|u^{k-2}-u^{k-3}\|^2.
\end{align}
Recall that in Lemma \ref{main-lem2}, we have proved
\begin{equation}
\sum_{k = 1}^n  d_{n-k+1} \left< L_{k}^\alpha u, D_k u\right> \geq \frac {5\alpha \Delta t^{1-\alpha}} {12\Gamma(3-\alpha)} \sum_{k=1}^{n}  d_{n-k+1} \left\| D_k u\right\|^2.
\end{equation}
We then have
\begin{equation}
\begin{aligned}
\sum_{k = 1}^n  d_{n-k+1} D_k E 
&\leq -\Delta t \sum_{k=1}^{n}  d_{n-k+1} \left(\frac {5\alpha } {12\Gamma(3-\alpha)\Delta t^{\alpha}} - 2 L\right)\left\| D_k u\right\|^2 \\
& + \Delta t \sum_{k=1}^{n-1}  d_{n-k}  L \left\| D_{k} u\right\|^2 
+ \Delta t \sum_{k=1}^{n-2}  d_{n-k-1} \frac L 2\left\| D_{k} u\right\|^2.
\end{aligned}
\end{equation}
Note that $4d_{j+1}>d_j$ according to Lemma \ref{lem1}.
When 
\begin{equation}
\Delta t^\alpha\leq \frac{5\alpha}{168\Gamma(3-\alpha)L},
\end{equation}
we then have 
\begin{equation}
\sum_{k = 1}^n  d_{n-k+1} D_k E \leq 0.
\end{equation}
\end{proof}

Theorem \ref{main-thm} gives a time-fractional energy law, which yields directly the following energy boundedness result for the L2 scheme \eqref{eq:L2AB} due to the decrease of $d_j$:
\begin{corollary}[Energy boundedness]\label{cor}
For the $3-\alpha$ order L2 scheme \eqref{eq:L2AB} with the same conditions in Theorem \ref{main-thm}, the following energy boundedness holds:
\begin{equation}
E^{n} \leq E^0, \quad \forall 1\leq n \leq N.
\end{equation}
\begin{proof}

	This theorem can be proved easily by mathematical induction. When n = 1, \eqref{main-thm-eq} is 
	\begin{align}
		d_1(E^1 - E^0)\leq0,
	\end{align}
	which indicates that $E^1\leq E^0$.
	Assuming that $ E^k\leq E^0$ for $1\leq k \leq n-1$ and rewritting \eqref{main-thm-eq} as
\begin{align}
	d_1 E^n \leq d_nE^0+\sum_{k=1}^{n-1}(-d_{n-k+1}+d_{n-k}) E^k.
\end{align}
Recalling that $d_j>0$ decreases w.r.t $j$, we have $E^n  \leq E^0.$
\end{proof}
\end{corollary}
\begin{remark}
One can also prove Corollary \ref{cor} directly using Lemma \ref{main-lem} and obtain a better restriction on $\Delta t$. We leave this prove to readers.
\end{remark}

\section{Numerical tests}
\label{sect4}
In this section, we test the proposed L2 schemes for time-fractional phase-field models, in order to verify the convergence rate and the energy stability.
More specifically, we consider the AC model with $\mathcal G = -1$ and the CH model with $\mathcal G = \Delta$.
The energy of the Allen--Cahn and Cahn--Hilliard equations is 
\begin{equation}\label{eq:energy}
E( u) = \int_\Omega \left(\frac{\varepsilon^2} 2 \left| \nabla  u \right|^2 + F( u) \right) \, {\rm d} x,
\end{equation}
where
\begin{equation}
F( u)  = \frac 1 4 \left(1- u^2\right)^2.
\end{equation}

\begin{example} \label{exam1}
{\em Consider the 2D fractional Allen-Cahn equation
\begin{equation}\label{eq:ACpoly}
\partial_t^\alpha u = \varepsilon^2\Delta u + u-u^3 + s\quad\mbox{on } [-\pi,\pi]^2\times (0,T],
\end{equation}
with periodic boundary condition and the source term $s(x,y,t)$ s.t. the exact solution is
\begin{equation}
u(x,y,t) = 0.2 t^5 \sin(x)\cos(y).
\end{equation}
}
\end{example}

In this test, we use the Fourier spectral method with $128\times 128$ modes for spatial discretization.
This number is large enough so that the spatial approximation error is negligible. 
We take $\varepsilon = 0.1$.
The errors and convergence rates are given in Table \ref{tab1} and \ref{tab2} computed respectively by the L2 SAV scheme \eqref{eq:L2SAV} and the implicit-explicit L2 scheme \eqref{eq:L2AB}.
It can be observed that \eqref{eq:L2SAV} is approximately second order and \eqref{eq:L2AB} is $3-\alpha$ order, as expected.

However, we emphasize that the convergence rates can be reached  when the exact solution is regular enough w.r.t. time. 
If not, graded time mesh might be needed to preserve the correct convergence order, see for example \cite{hou2021highly} for some interesting discussions.

\begin{table}[htb!]
\renewcommand\arraystretch{1.4}
\begin{center}
\def\temptablewidth{1.1\textwidth}
\caption{$\ell_2$-errors at $T = 1$ for Example \ref{exam1} for $\alpha = 0.1$ (top) and $0.9$ (bottom) and their convergence rates, computed by the L2 SAV scheme.}\label{tab1}
{\rule{\temptablewidth}{1pt}}
\begin{tabular*}{\temptablewidth}{@{\extracolsep{\fill}}ccccccc}
   $\tau$     &$\frac1{40}$     &$\frac1{80}$
   &$\frac1{160}$  & $\frac1{320}$ & $\frac1{640}$ & $\frac1{1280}$ \\ \hline
  $\ell_2$-error   & $ 3.4147\times 10^{-2}$   & $8.9402\times 10^{-3}$   & $2.2826\times 10^{-3}$ & $5.7686\times 10^{-4}$ & $1.4502\times 10^{-4}$ & $3.6357\times 10^{-5}$ \\[3pt]
rate  & -- &$1.9334$ &$1.9696$ &$1.9844$ & $1.9920$ & $1.9959$
\end{tabular*}
{\rule{\temptablewidth}{1pt}}
{\rule{\temptablewidth}{1pt}}
\begin{tabular*}{\temptablewidth}{@{\extracolsep{\fill}}ccccccc}
   $\tau$     &$\frac1{40}$     &$\frac1{80}$
   &$\frac1{160}$  & $\frac1{320}$ & $\frac1{640}$ & $\frac1{1280}$ \\ \hline
     $\ell_2$-error   & $4.1677\times 10^{-4}$   & $1.6061\times 10^{-4}$   & $5.1724\times 10^{-5}$ & $1.5388\times 10^{-5}$ & $4.3863\times 10^{-6}$ & $1.2177\times 10^{-6}$ \\[3pt]
rate  & -- &$1.3757$ &$1.6346$ &$1.7490$ & $1.8108$ & $1.8488$
\end{tabular*}
{\rule{\temptablewidth}{1pt}}
\end{center}

\renewcommand\arraystretch{1.4}
\begin{center}
\def\temptablewidth{1.1\textwidth}
\caption{$\ell_2$-errors at $T = 1$ for Example \ref{exam1} for $\alpha = 0.1$ (top) and $0.9$ (bottom), and their convergence rates, computed by the implicit-explicit L2 scheme.}\label{tab2}
{\rule{\temptablewidth}{1pt}}
\begin{tabular*}{\temptablewidth}{@{\extracolsep{\fill}}ccccccc}
   $\tau$     &$\frac1{40}$     &$\frac1{80}$
   &$\frac1{160}$  & $\frac1{320}$ & $\frac1{640}$ & $\frac1{1280}$ \\ \hline
  $\ell_2$-error   & $2.1833\times 10^{-3}$   & $2.6112\times 10^{-4}$   & $3.1957\times 10^{-5}$ & $3.9571\times 10^{-6}$ & $4.9335\times 10^{-7}$ & $6.1750\times 10^{-8}$ \\[3pt]
rate  & -- &$3.0637$ &$3.0305$ &$3.0136$ & $3.0037$ & $2.9981$
\end{tabular*}
{\rule{\temptablewidth}{1pt}}
{\rule{\temptablewidth}{1pt}}
\begin{tabular*}{\temptablewidth}{@{\extracolsep{\fill}}ccccccc}
   $\tau$     &$\frac1{40}$     &$\frac1{80}$
   &$\frac1{160}$  & $\frac1{320}$ & $\frac1{640}$ & $\frac1{1280}$ \\ \hline
     $\ell_2$-error   & $1.9656\times 10^{-3}$   & $4.9721\times 10^{-4}$   & $1.2088\times 10^{-4}$ & $2.8802\times 10^{-5}$ & $6.7924\times 10^{-6}$ & $1.5934\times 10^{-6}$ \\[3pt]
rate  & -- &$1.9830$ &$2.0402$ &$2.0694$ & $2.0842$ & $2.0918$
\end{tabular*}
{\rule{\temptablewidth}{1pt}}
\end{center}
\end{table}

\begin{example} \label{exam2}
{\em Consider the 2D fractional Allen-Cahn equation
\begin{equation}\label{eq:ACpoly}
\partial_t^\alpha u = \varepsilon^2\Delta u + u-u^3 \quad\mbox{on } [0,2\pi]^2\times (0,T],
\end{equation}
with periodic boundary condition and initial condition composed of seven circles with centers and radii given in Table \ref{tab:xyr}:
\begin{equation}\label{eq:init_7circ}
u_0(x,y) = -1 + \sum_{i=1}^7 f\left( \sqrt{(x-x_i)^2+(y-y_i)^2} - r_i\right),
\end{equation}
where
\begin{equation}
f(s) = \left\{
\begin{aligned}
& 2 e^{-\varepsilon^2/s^2} &&\mbox{if } s<0,\\
& 0 && \mbox{otherwise.}
\end{aligned}
\right.
\end{equation}
}
\end{example}

\begin{table}[htb!]
\renewcommand\arraystretch{1.3}
\begin{center}
\def\temptablewidth{0.8\textwidth}
\caption{Centers $(x_i,y_i)$ and radii $r_i$ in the initial condition \eqref{eq:init_7circ}, which are the same as in \cite{church2019high}.}\label{tab:xyr}
{\rule{\temptablewidth}{1.1pt}}
\begin{tabular*}{\temptablewidth}{@{\extracolsep{\fill}}c|ccccccc}
   $i$   &1 & 2     &3   &4  & 5 & 6 & 7 \\  \hline
  $x_i$  & $\pi/2$   & $\pi/4$   & $\pi/2$ & $\pi$ & $3\pi/2$ & $\pi$ & $3\pi/2$ \\[3pt]
$y_i$  & $\pi/2$   & $3\pi/4$   & $5\pi/4$ & $\pi/4$ & $\pi/4$ & $\pi$ & $3\pi/2$ \\[3pt]
$r_i$  & $\pi/5$   & $2\pi/15$   & $2\pi/15$ & $\pi/10$ & $\pi/10$ & $\pi/4$ & $\pi/4$ \\[3pt]
\end{tabular*}
{\rule{\temptablewidth}{1.1pt}}
\end{center}
\end{table}

\begin{figure}[!]
\centering
\includegraphics[trim = {0in 0.8in 0.in 0},clip,width = 0.99\textwidth]{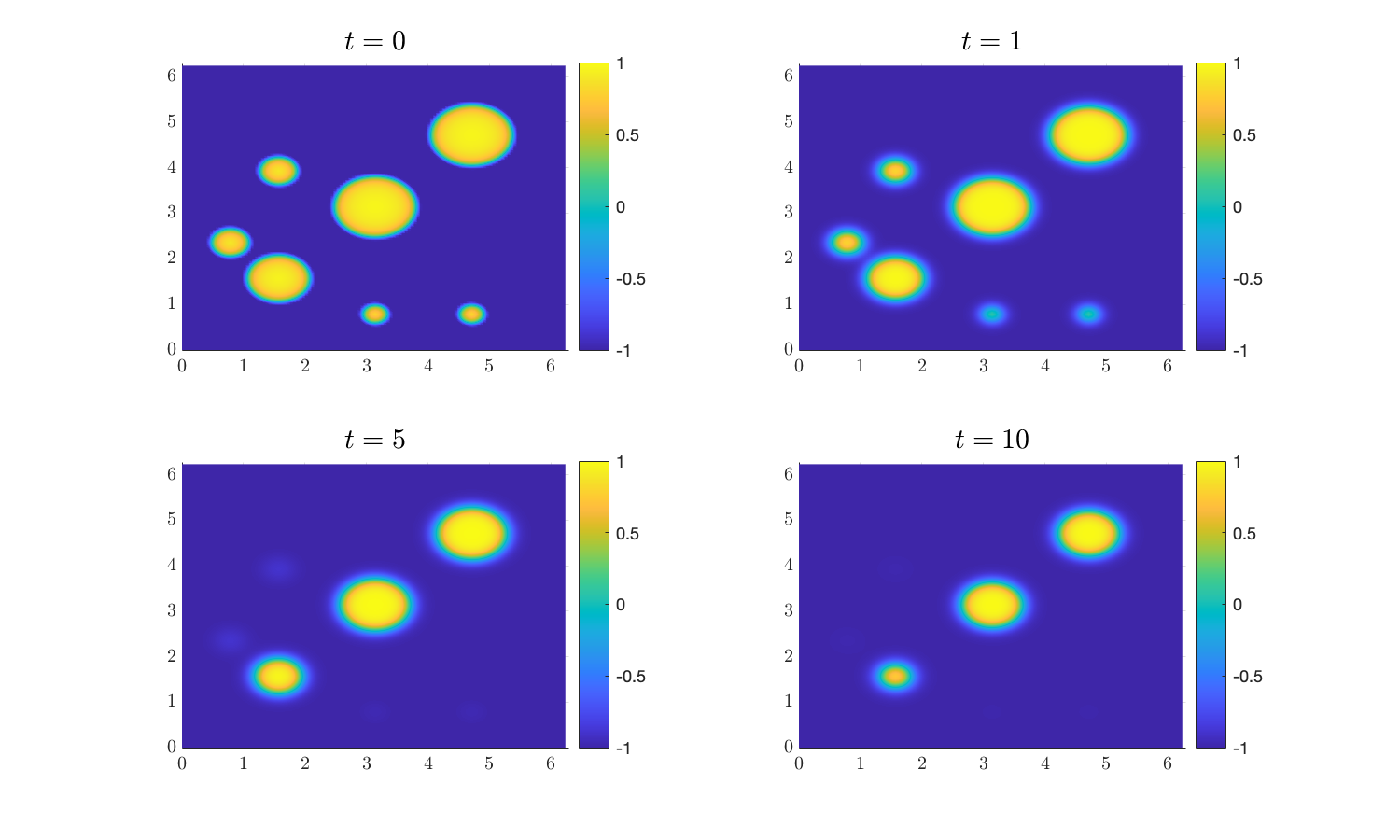}
\caption{Numerical solution of Example \ref{exam2} with $\alpha = 0.9, ~\Delta t = 0.01$ and the number of Fourier modes $128\times 128$, computed by the implicit-explicit L2 scheme.}\label{fig:sol_7circ}
\includegraphics[trim = {0in 0.in 0.5in 0.2in},clip,width = 0.5\textwidth]{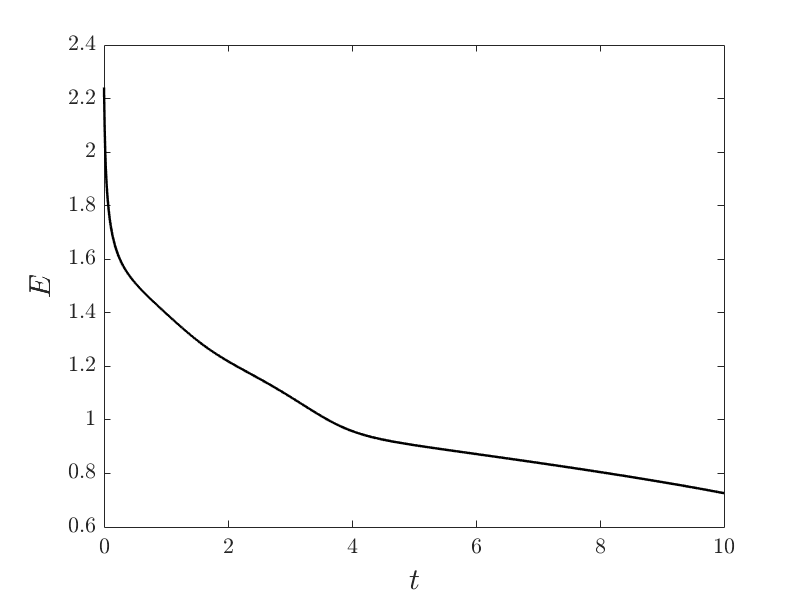}
\caption{Classical energy w.r.t. time for Example \ref{exam2}  with $\alpha = 0.9, ~\Delta t = 0.01$ and the number of Fourier modes $128\times 128$, computed by the implicit-explicit L2 scheme.}\label{fig:energy_7circs}
\end{figure}

We take $\varepsilon = 0.1,~\alpha = 0.9,~\Delta t = 0.01$ and use $128\times 128$ Fourier modes for spatial discretization.
The numerical solution and energy evolution are illustrated respectively in Figure \ref{fig:sol_7circ} and \ref{fig:energy_7circs}.
In this case, we can observe that the classical energy decreases w.r.t. time.

\begin{example} \label{exam3}
{\em Consider the 2D fractional Cahn--Hilliard equation
\begin{equation}\label{eq:CHpoly}
\partial_t^\alpha u = -\varepsilon^2\Delta^2 u + \Delta(u-u^3) \quad\mbox{on } [0,2\pi]^2\times (0,T],
\end{equation}
with periodic boundary condition and random initial condition distributed uniformly in $[- 0.5, 0.5]$.
}
\end{example}

\begin{figure}[!]
\centering
\includegraphics[trim = {0in 0.8in 0.in 0},clip,width = 0.99\textwidth]{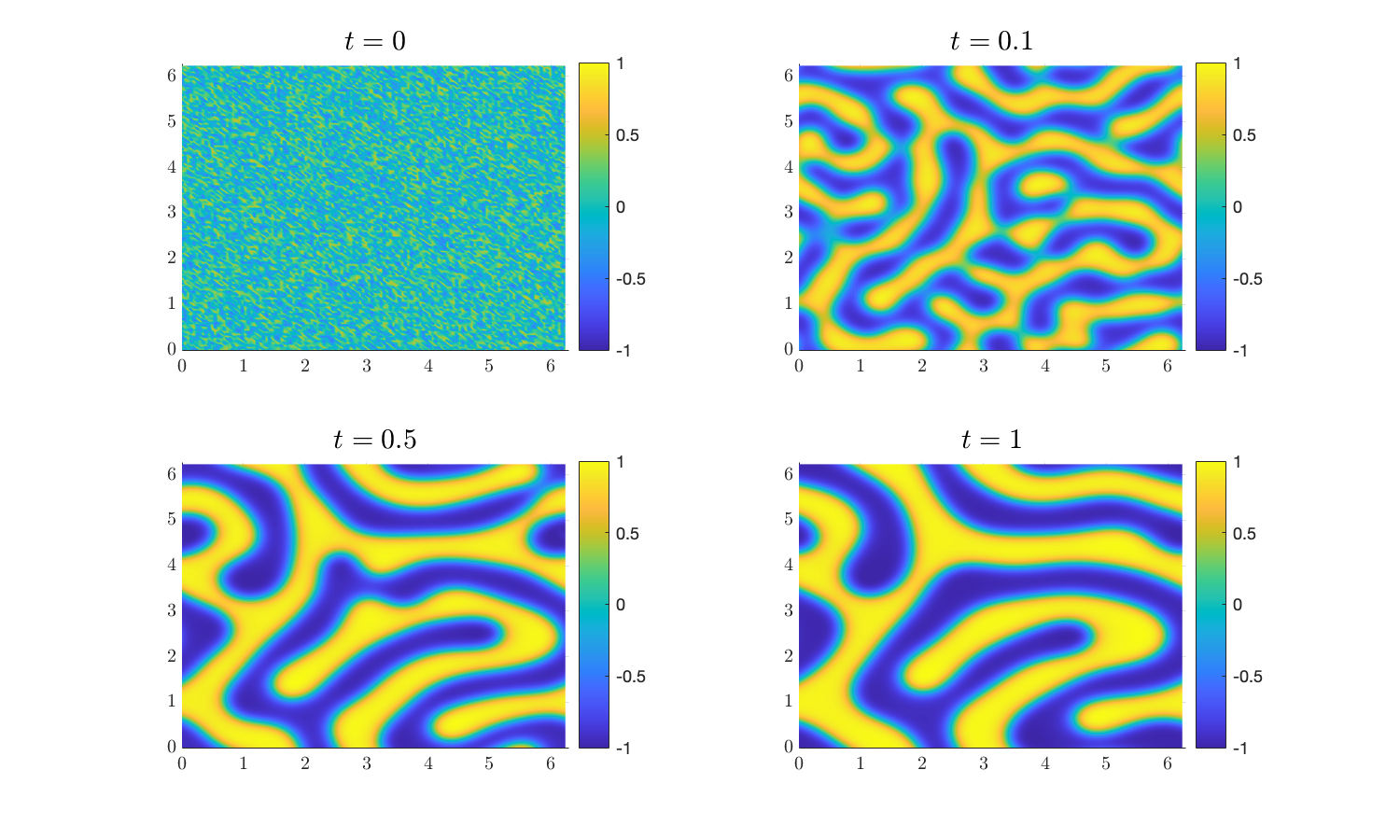}
\caption{Numerical solution of Example \ref{exam2} with $\alpha = 0.8, ~\Delta t = 0.001$ and the number of Fourier modes $128\times 128$, computed by the implicit-explicit L2 scheme.}\label{fig:ch_rand}
\includegraphics[trim = {0in 0.in 0.5in 0.2in},clip,width = 0.5\textwidth]{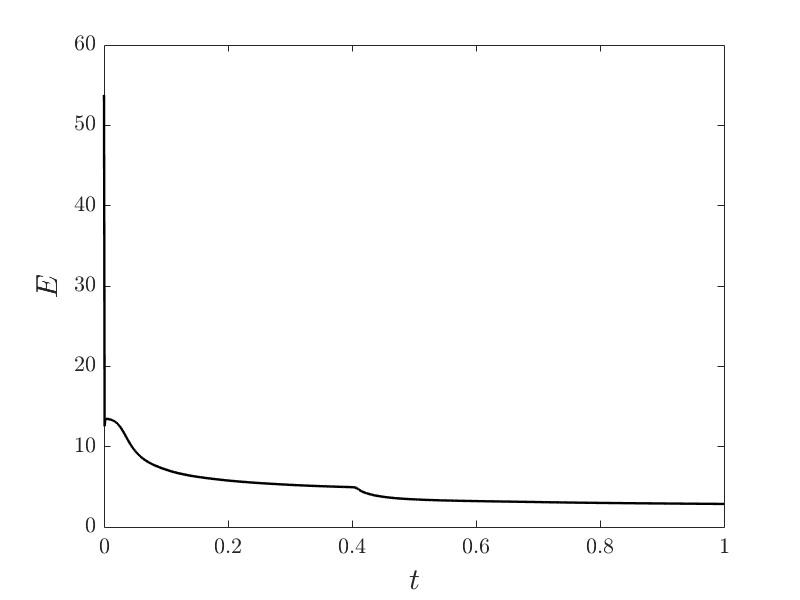}
\caption{Classical energy w.r.t. time for Example \ref{exam2}  with $\alpha = 0.8, ~\Delta t = 0.001$ and the number of Fourier modes $128\times 128$, computed by the implicit-explicit L2 scheme.}\label{fig:energy_ch_rand}
\end{figure}

We take $\varepsilon = 0.1,~\alpha = 0.8,~\Delta t = 0.001,~T = 1$, and use $128\times 128$ Fourier modes for spatial discretization.
The numerical solution and energy evolution are illustrated respectively in Figure \ref{fig:ch_rand} and \ref{fig:energy_ch_rand}.
It can be observed that the modified energy is bounded by initial energy. 
Note that near $t=0$, the energy dissipation property seems destroyed but the energy boundedness is still satisfied.
Similar situation has also been reported in \cite{hou2019variant}.

\section{Conclusion}
We have established the energy boundedness of the second order L2 SAV scheme for any phase-field equation and the time-fractional energy law of the $3-\alpha$ order L2 IMEX scheme for the AC equation.
To prove the energy stability, a reformulation of L2 approximation is proposed and several useful properties have been provided for the L2 operator.
Numerical tests are provided to verify the convergence order (when the exact solution is sufficiently regular w.r.t. time) and the energy stability.  

However, we shall mention that it is still an open question whether the rigorous energy dissipation holds (even on the continuous level), which is challenging due to the existence of both nonlocality and nonlinearity. 

\section*{Acknowledgements}
%This work is partially supported by the Special Project on High-Performance Computing of the National Key R\&D Program under No. 2016YFB0200604, the National Natural Science Foundation of China (NSFC) Grant No. 11731006, the NSFC/Hong Kong RRC Joint Research Scheme (NSFC/RGC 11961160718), and the fund of the Guangdong Provincial Key Laboratory of Computational Science and Material Design (No. 2019B030301001). The work of J. Yang is supported by the National Science Foundation of China (NSFC-11871264) and the Guangdong Basic and Applied Basic Research Foundation (2018A0303130123). 
The research of C. Quan is supported by NSFC Grant 11901281, the Guangdong Basic and Applied Basic Research Foundation (2020A1515010336), and the Stable Support Plan Program of Shenzhen Natural Science Fund (Program Contract No. 20200925160747003).

\appendix
\renewcommand{\theequation}{\thesection.\arabic{equation}}
\section{Signs of $Q_1$, $Q_2$, $Q_3$ in \eqref{eq:Q}}\label{append}
For the simplicity, we denote 
\begin{equation}
\theta_1 = \frac{1}{i} \quad \mbox{and} \quad \theta_2 = \frac{1}{i-j+1},
\end{equation}
so that $0<\theta_1 \leq \theta_2\leq \frac 1 2$ since $i>j\geq 1$.
Then, we can rewrite $\kappa$ define in \eqref{eq:kappa}as 
\begin{equation}\label{eq:kappa12}
  \begin{array}{r@{}l}
	\begin{aligned}
\kappa(i,\beta) & = i^\beta\rho(\theta_1,\beta),\\
\kappa(i-j+1,\beta) & = (i-j+1)^\beta\rho(\theta_2,\beta).
	 \end{aligned}
  \end{array}
\end{equation}
with
\begin{equation}\label{eq:Taylor}
\rho(\theta,\beta) \coloneqq (1+\theta)^\beta - 2 + (1-\theta)^\beta = 2 \sum_{m=1}^\infty \binom{\beta}{2m} \theta^{2m}.
%+ 2 \binom{\beta}{4} \theta^4 + 2\binom{\beta}{6}\theta^6 + \cdots.
\end{equation}

Firstly, we prove that $Q_1\geq 0$ in \eqref{eq:Q}. 
Combining the first equation of \eqref{eq:Q}, \eqref{eq:kappa12}, and \eqref{eq:Taylor}, we have
\begin{equation}
  \begin{array}{r@{}l}
	\begin{aligned}
	Q_1 & =  \frac \alpha 2 \kappa(i,2-\alpha) \kappa(i-j+1,-\alpha-1) + \kappa(i,-\alpha) \kappa(i-j+1,1-\alpha)\\
%	& = 2 \alpha \,  i^{2-\alpha} (i-j+1)^{-\alpha-1} \sum_{m=1}^\infty \binom{2-\alpha}{2m} \theta_1^{2m} \sum_{m=1}^\infty \binom{-\alpha-1}{2m} \theta_2^{2m} \\
%	& + 4 i^{-\alpha} (i-j+1)^{1-\alpha} \sum_{m=1}^\infty \binom{-\alpha}{2m} \theta_1^{2m} \sum_{m=1}^\infty \binom{1-\alpha}{2m} \theta_2^{2m} \\
%	& \geq 4  i^{-\alpha} (i-j+1)^{1-\alpha} \Big [ - \sum_{m=1}^\infty \frac{(2-\alpha)}{(\alpha+2m-2)(\alpha+2m-1)}\binom{-\alpha}{2m} \theta_1^{2m} \\
%	& \sum_{m=1}^\infty \frac{(\alpha+2m-1)(\alpha+2m)}{2}\binom{1-\alpha}{2m} \theta_2^{2m} 
%	 + \sum_{m=1}^\infty \binom{-\alpha}{2m} \theta_1^{2m} \sum_{m=1}^\infty \binom{1-\alpha}{2m} \theta_2^{2m}\Big] \\
	& \geq \frac 1 2  i^{-\alpha} (i-j+1)^{1-\alpha} \Big [ \alpha\rho(\theta_1,2-\alpha)\rho(\theta_2,-\alpha-1) + 2\rho(\theta_1,-\alpha)\rho(\theta_2,1-\alpha)\Big]\\
	& = \frac 1 2  i^{-\alpha} (i-j+1)^{1-\alpha} H_1,
	\end{aligned}
  \end{array}
\end{equation}
with
\begin{equation}
  \begin{array}{r@{}l}
	\begin{aligned}
	H_1& = \alpha\rho(\theta_1,2-\alpha)\rho(\theta_2,-\alpha-1) + 2\rho(\theta_1,-\alpha)\rho(\theta_2,1-\alpha).
	\end{aligned}
  \end{array}
\end{equation} 
As $\theta_2\leq \frac 1 2$, it is not difficult to verify
\begin{equation}\label{eq:half}
(\alpha+1) \geq  \sum_{m = 2}^\infty (2m-3)(\alpha+1)\theta_2^{2m-2} \geq \sum_{m=2}^\infty  \frac{2(2m-3)}{2m+1}\binom{-\alpha-1}{2m} \theta_2^{2m-2},
\end{equation}
due to the fact that
\begin{equation}
(\alpha+1)\geq \frac{2}{2m+1} \binom{-\alpha-1}{2m} =\frac{2(\alpha+1)\cdots(\alpha+2m)}{(2m+1)!} .
\end{equation}
Combining \eqref{eq:Taylor} and \eqref{eq:half}, we derive
\begin{equation}
  \begin{array}{r@{}l}
	\begin{aligned}
	\rho(\theta_2,-\alpha-1) 
	& =  2 \sum_{m=1}^\infty \binom{-\alpha-1}{2m} \theta_2^{2m}  = (\alpha+1)(\alpha+2) \theta_2^2 + \sum_{m=2}^\infty 2\binom{-\alpha-1}{2m} \theta_2^{2m}  \\
	& \geq (\alpha+1)^2 \theta_2^2  + 4 \sum_{m=2}^\infty  \frac{2m-1}{2m+1}\binom{-\alpha-1}{2m} \theta_2^{2m}.
	\end{aligned}
  \end{array}
\end{equation}
As a consequence, we have
\begin{equation}\label{eq:htheta1}
  \begin{array}{r@{}l}
	\begin{aligned}
H_1 &
\geq 2 \alpha \sum_{m_1=1}^\infty \binom{2-\alpha}{2m_1} \theta_1^{2m_1} \left[(\alpha+1)^2 \theta_2^2  + 4 \sum_{m_2=2}^\infty  \frac{2m_2-1}{2m_2+1}\binom{-\alpha-1}{2m_2} \theta_2^{2m_2}\right]\\
&+ 8\sum_{m_1=1}^\infty \binom{-\alpha}{2m_1} \theta_1^{2m_1} \sum_{m_2=1}^\infty \binom{1-\alpha}{2m_2} \theta_2^{2m_2} \\
& = 8 \sum_{m_1=1}^\infty \sum_{m_2=1}^\infty c_{m_1,m_2}\theta_1^{2m_1}\theta_2^{2m_2}.
	\end{aligned}
  \end{array}
\end{equation}
In the case of $m_1 = m_2= m$, we can find  that if $m = 1$,
\begin{equation}
  \begin{array}{r@{}l}
	\begin{aligned}
	c_{1,1} = \frac\alpha 4 (\alpha+1)^2\binom{2-\alpha}{2}  + \binom{-\alpha}{2}\binom{1-\alpha}{2} = \frac{\alpha(1-\alpha)^2(\alpha+1)(\alpha+2)}{8}\geq 0,
	\end{aligned}
  \end{array}
\end{equation}
while if $ m\geq 2$,
\begin{equation}
  \begin{array}{r@{}l}
	\begin{aligned}
	c_{m,m} & =   \frac{\alpha(2m-1)}{2m+1}\binom{2-\alpha}{2m} \binom{-\alpha-1}{2m}+ \binom{-\alpha}{2m}\binom{1-\alpha}{2m} \\
	& = \binom{2-\alpha}{2m} \binom{-\alpha-1}{2m} \left[ \frac{\alpha(2m-1)}{2m+1} -\frac{\alpha(2m-2+\alpha)}{(2-\alpha)(2m+\alpha)}\right] \\
	& = \binom{2-\alpha}{2m} \binom{-\alpha-1}{2m} \alpha \left[ \frac{1-\alpha}{2-\alpha}+\frac 2 {2m+\alpha} -\frac 2 {2m+1}\right] \\
	&\geq 0.
	\end{aligned}
  \end{array}
\end{equation}
In the case of $m_1> m_2=1$, we have
\begin{equation}
  \begin{array}{r@{}l}
	\begin{aligned}
	& c_{m_1,1} \theta_1^{2m_1} \theta_2^{2} + c_{1,m_1} \theta_1^{2} \theta_2^{2m_1} \\
	& = \left[\frac\alpha 4 (\alpha+1)^2\binom{2-\alpha}{2m_1}  + \binom{-\alpha}{2m_1}\binom{1-\alpha}{2}\right] \theta_1^{2m_1}\theta_2^2 \\
	& + \left[ \frac{\alpha(2m_1-1)}{2m_1+1}\binom{2-\alpha}{2} \binom{-\alpha-1}{2m_1}+ \binom{-\alpha}{2}\binom{1-\alpha}{2m_1} \right] \theta_1^{2} \theta_2^{2m_1} \\
	& \geq \frac \alpha 4 \bigg[(\alpha+1)^2 - 2(2m_1+\alpha-1)(2m_1+\alpha-2) + \frac{2(2m_1-1)}{2m_1+1} \left(\frac{2m_1}{\alpha}+1\right) \\
	& \quad (2m_1+\alpha-1)(2m_1+\alpha-2) 
	-2(2m_1+\alpha-2)\bigg] \binom{2-\alpha}{2m_1}\theta_1^{2m_1}\theta_2^2 \\
	& \geq 0,
	\end{aligned}
  \end{array}
\end{equation}
where we use the fact $\theta_1\leq \theta_2$ and $m_1\geq m_2+1 = 2$.
In the case of $m_1> m_2\geq 2$, we have
\begin{align*}\label{eq:cm1m2}
	& c_{m_1,m_2} \theta_1^{2m_1} \theta_2^{2m_2} + c_{m_2,m_1} \theta_1^{2m_2} \theta_2^{2m_1} \\
	& =\left[ \frac{\alpha(2m_2-1)}{2m_2+1}\binom{2-\alpha}{2m_1} \binom{-\alpha-1}{2m_2}+ \binom{-\alpha}{2m_1}\binom{1-\alpha}{2m_2} \right] \theta_1^{2m_1} \theta_2^{2m_2}  \\
	& + \left[ \frac{\alpha(2m_1-1)}{2m_1+1}\binom{2-\alpha}{2m_2} \binom{-\alpha-1}{2m_1}+ \binom{-\alpha}{2m_2}\binom{1-\alpha}{2m_1} \right] \theta_1^{2m_2} \theta_2^{2m_1} \\
	& =\alpha\left[ \frac{2m_2-1}{2m_2+1}-\frac{(2m_1-1+\alpha)(2m_1-2+\alpha)}{(2-\alpha)(2m_2+\alpha)(2m_2-1+\alpha)}\right] \binom{2-\alpha}{2m_1} \binom{-\alpha-1}{2m_2} \theta_1^{2m_1} \theta_2^{2m_2}  \\
	& + \alpha\left[ \frac{2m_1-1}{2m_1+1}-\frac{(2m_2-1+\alpha)(2m_2-2+\alpha)}{(2-\alpha)(2m_1+\alpha)(2m_1-1+\alpha)}\right] \binom{2-\alpha}{2m_2} \binom{-\alpha-1}{2m_1} \theta_1^{2m_2} \theta_2^{2m_1} \\
	&\geq \alpha \binom{2-\alpha}{2m_1} \binom{-\alpha-1}{2m_2} \theta_1^{2m_1} \theta_2^{2m_2} \bigg[ \frac{2m_2-1}{2m_2+1}-\frac{(2m_1-1+\alpha)(2m_1-2+\alpha)}{(2-\alpha)(2m_2+\alpha)(2m_2-1+\alpha)} \\
	& + \frac{(2m_1-1)(2m_1+\alpha)(2m_1-1+\alpha)(2m_1-2+\alpha)}{(2m_1+1)(2m_2+\alpha)(2m_2-1+\alpha)(2m_2-2+\alpha)}
	- \frac{2m_1-2+\alpha}{(2-\alpha)(2m_2+\alpha)} \bigg] \\
	& \geq \alpha \binom{2-\alpha}{2m_1} \binom{-\alpha-1}{2m_2} \theta_1^{2m_1} \theta_2^{2m_2} \bigg[\frac{2m_2-1}{2m_2+1} - \frac{2m_1-2+\alpha}{(2m_2+\alpha)} \\
	& +\frac{(2m_1-1+\alpha)(2m_1-2+\alpha)}{(2m_2+\alpha)(2m_2-1+\alpha)}\left( \frac{(2m_1-1)(2m_1+\alpha)}{(2m_1+1)(2m_2-2+\alpha)} -1\right)\bigg] \\
	& \geq \alpha \binom{2-\alpha}{2m_1} \binom{-\alpha-1}{2m_2} \theta_1^{2m_1} \theta_2^{2m_2} \bigg[- \frac{2(m_1-m_2)}{(2m_2+\alpha)} + \frac{(2m_1-1)(2m_1+\alpha)}{(2m_1+1)(2m_2-2+\alpha)} -1\bigg]  \\
	%& \geq \alpha \binom{2-\alpha}{2m_1} \binom{-\alpha-1}{2m_2} \theta_1^{2m_1} \theta_2^{2m_2} \bigg[- \frac{2(m_1-m_2)}{(2m_2+\alpha)} + \frac{2(1-\alpha)}{(2m_1+1)(2m_2-2+\alpha)} + \frac{2(m_1-m_2)}{(2m_2-2+\alpha)}\bigg] \\
	& \geq 0. \numberthis
\end{align*}
Combining \eqref{eq:htheta1}--\eqref{eq:cm1m2}, we then claim
$H_1\geq 0$, which yields $Q_1 \geq 0$.

Secondly, we prove that $Q_2\leq 0$ in \eqref{eq:Q}. 
Combining the second equation of \eqref{eq:Q}, \eqref{eq:kappa12}, and \eqref{eq:Taylor}, we have
\begin{equation}
  \begin{array}{r@{}l}
	\begin{aligned}
	Q_2 & = \frac \alpha 2 \kappa(i,1-\alpha) \kappa(i-j+1,-\alpha-1)+ \frac 1 2 (1-\alpha) \kappa(i,-\alpha) \kappa(i-j+1,-\alpha),\\
	& \leq  \frac 1 2  i^{1-\alpha} (i-j+1)^{-\alpha-1} \Big [ \alpha\rho(\theta_1,1-\alpha)\rho(\theta_2,-\alpha-1) + (1-\alpha)\rho(\theta_1,-\alpha)\rho(\theta_2,-\alpha)\Big]\\
	& = \frac 1 2  i^{1-\alpha} (i-j+1)^{-\alpha-1}  H_2,
	\end{aligned}
  \end{array}
\end{equation}
with 
\begin{equation}
  \begin{array}{r@{}l}
	\begin{aligned}
H_2 & = \alpha\rho(\theta_1,1-\alpha)\rho(\theta_2,-\alpha-1) + (1-\alpha)\rho(\theta_1,-\alpha)\rho(\theta_2,-\alpha) \\
& = 4 \alpha \sum_{m_1=1}^\infty \binom{1-\alpha}{2m_1} \theta_1^{2m_1}  \sum_{m_2=1}^\infty \binom{-\alpha-1}{2m_2} \theta_2^{2m_2}+ 4 (1-\alpha)\sum_{m_1=1}^\infty \binom{-\alpha}{2m_1} \theta_1^{2m_1} \sum_{m_2=1}^\infty \binom{-\alpha}{2m_2} \theta_2^{2m_2} \\
& = 4 \alpha^2 (1-\alpha)(1+\alpha)\bigg[ - \sum_{m_1=1}^\infty \frac{1}{2m_1(2m_1-1)}\binom{-\alpha-1}{2m_1-2} \theta_1^{2m_1}
\sum_{m_2=1}^\infty  \frac{\alpha+2}{2m_2(2m_2-1)}\binom{-\alpha-3}{2m_2-2} \theta_2^{2m_2} \\
&+ \sum_{m_1=1}^\infty \frac{1}{2m_1(2m_1-1)}\binom{-\alpha-2}{2m_1-2} \theta_1^{2m_1}
\sum_{m_2=1}^\infty \frac{\alpha+1}{2m_2(2m_2-1)}\binom{-\alpha-2}{2m_2-2} \theta_2^{2m_2}\bigg]\\
& \leq 0,
	\end{aligned}
  \end{array}
\end{equation}
where the proof of the last inequality is similar to the previous case of $H_2$.
As a consequence, we can claim that $Q_2\leq 0$.

Thirdly, we prove that $Q_3\leq 0$ in \eqref{eq:Q}. 
Combining the second equation of \eqref{eq:Q}, \eqref{eq:kappa12}, and \eqref{eq:Taylor}, we have
\begin{equation}
  \begin{array}{r@{}l}
	\begin{aligned}
	Q_3 & = -(1-\alpha) \kappa(i,2-\alpha) \kappa(i-j+1,-\alpha) + (2-\alpha) \kappa(i,1-\alpha) \kappa(i-j+1,1-\alpha),\\
	& \leq  \frac 1 2  i^{2-\alpha} (i-j+1)^{-\alpha} \Big [ -(1-\alpha)\rho(\theta_1,2-\alpha)\rho(\theta_2,-\alpha) + (2-\alpha)\rho(\theta_1,1-\alpha)\rho(\theta_2,1-\alpha)\Big]\\
	& = \frac 1 2  i^{2-\alpha} (i-j+1)^{-\alpha}  H_3,
	\end{aligned}
  \end{array}
\end{equation}
with 
\begin{equation}
  \begin{array}{r@{}l}
	\begin{aligned}
H_3 & =  -(1-\alpha)\rho(\theta_1,2-\alpha)\rho(\theta_2,-\alpha) + (2-\alpha)\rho(\theta_1,1-\alpha)\rho(\theta_2,1-\alpha)\\
& = -4 (1-\alpha) \sum_{m_1=1}^\infty \binom{2-\alpha}{2m_1} \theta_1^{2m_1}  \sum_{m_2=1}^\infty \binom{-\alpha}{2m_2} \theta_2^{2m_2}+ 4 (2-\alpha)\sum_{m_1=1}^\infty \binom{1-\alpha}{2m_1} \theta_1^{2m_1} \sum_{m_2=1}^\infty \binom{1-\alpha}{2m_2} \theta_2^{2m_2} \\
& = 4 \alpha (1-\alpha)^2 (2-\alpha)\bigg[-\sum_{m_1=1}^\infty \frac{1}{2m_1(2m_1-1)}\binom{-\alpha}{2m_1-2} \theta_1^{2m_1}
\sum_{m_2=1}^\infty \frac{\alpha+1}{2m_2(2m_2-1)}\binom{-\alpha-2}{2m_2-2} \theta_2^{2m_2} \\
& + \sum_{m_1=1}^\infty \frac{1}{2m_1(2m_1-1)}\binom{-\alpha-2}{2m_1-2} \theta_1^{2m_1}
\sum_{m_2=1}^\infty \frac{\alpha}{2m_2(2m_2-1)}\binom{-\alpha-1}{2m_2-2} \theta_2^{2m_2}\bigg]\\
& \leq 0,
	\end{aligned}
  \end{array}
\end{equation}
where the proof of the last inequality is similar to the previous case of $H_2$.
As a consequence, we can claim that $Q_3\leq 0$.

\bibliography{bibfile}
\bibliographystyle{unsrt}

\end{document}